\documentclass[preprint]{elsarticle} 
\usepackage{lineno,hyperref}
\modulolinenumbers[5] 

\journal{Acta Astronautica}

\usepackage{bm}
\usepackage{amsmath}
\usepackage{subfigure}
\usepackage{overcite}
\usepackage{footnpag}
\usepackage[justification=justified]{caption}
\usepackage{wasysym}
\usepackage{amssymb}
\usepackage{array}
\usepackage{multiro w}
\usepackage{epstopdf}
\usepackage{graphicx}
\usepackage{soul}
\usepackage{gensymb}

\usepackage{multirow}
\usepackage{booktabs,caption}
\usepackage[flushleft]{threeparttable}
\usepackage{hhline}
\usepackage[font=small,skip=0pt]{caption}
\usepackage[space]{grffile}
\usepackage{marvosym}
\usepackage{chngpage}

\usepackage{comment}
\usepackage{float}

\usepackage{lipsum}
\usepackage{booktabs}
\usepackage{adjustbox}

\usepackage{tabularx}
\usepackage{multirow}

\usepackage{siunitx}
\usepackage{graphicx,subfigure}
\usepackage{ragged2e}
\usepackage{xcolor}
\usepackage{mathtools}

\begin{document}

\begin{frontmatter}

\title{Analytical Process Noise Covariance Modeling for \protect\\ Absolute and Relative Orbits} 

\author[]{Nathan Stacey\corref{mycorrespondingauthor}}
\cortext[mycorrespondingauthor]{Corresponding author}
\ead{nstacey@stanford.edu}
\author[]{Simone D'Amico}
\address{Department of Aeronautics and Astronautics, Stanford University,
Durand Building, 496 Lomita Mall, Stanford, CA 94305-4035, USA}

\begin{abstract}
This paper develops new analytical process noise covariance models for both absolute and relative spacecraft states. Process noise is always present when propagating a spacecraft state due to dynamics modeling deficiencies. Accurately modeling this noise is essential for sequential orbit determination and improves satellite conjunction analysis. A common approach called state noise compensation models process noise as zero-mean Gaussian white noise accelerations. The resulting process noise covariance can be evaluated numerically, which is computationally intensive, or through a widely used analytical model that is restricted to an absolute Cartesian state and small propagation intervals. Moreover, mathematically rigorous, analytical process noise covariance models for relative spacecraft states are not currently available. To address these limitations of the state of the art, new analytical process noise covariance models are developed for state noise compensation for both Cartesian and orbital element absolute state representations by leveraging spacecraft relative dynamics models. Two frameworks are then presented for modeling the process noise covariance of relative spacecraft states by assuming either small or large interspacecraft separation. The presented techniques are validated through numerical simulations.

\end{abstract}


\end{frontmatter}


\section{Introduction}
In Kalman filtering, process noise represents deviations of the true dynamics from the modeled dynamics. Process noise is always present in orbit determination due to imperfect knowledge of the numerous complex forces acting on a spacecraft and limited computational resources for modeling those forces. Computational capabilities are especially constrained for onboard orbit determination, often necessitating reduced order dynamics models. Accurately capturing dynamics modeling deficiencies through modeled process noise is essential for optimal sequential estimation. Inaccurate process noise modeling increases estimation error and can result in filter inconsistency and divergence\cite{schlee_divergence_1967,mehra_approaches_1972}. Realistic process noise modeling improves the fidelity of the spacecraft state error covariance, which increases the accuracy of satellite conjunction analysis. In particular, analytical process noise models are advantageous because they are computationally efficient and provide insight into system behavior. This paper addresses gaps in the state of the art through the development of new analytical process noise covariance models for absolute and relative spacecraft states parameterized using both Cartesian coordinates and orbital elements. Notably, analytical process noise covariance models for absolute spacecraft states are obtained by exploiting relative spacecraft dynamics models to describe the state error dynamics.

A variety of approaches have been established to account for unmodeled accelerations in orbit determination. Here the term unmodeled accelerations refers to the difference between the true and modeled spacecraft accelerations, which is due to forces that are omitted in the modeled dynamics as well as errors in the modeled forces. As demonstrated by the TOPEX/Poseidon and Jason-1 missions, unmodeled accelerations can be approximated using trigonometric functions with fixed frequencies and estimated coefficients\cite{luthcke_1-centimeter_2003,marshall_temporal_1995}. However, this technique is best suited to a posteriori high-precision orbit determination\cite{montenbruck_satellite_force_2000}. In batch estimation, it is also possible to estimate piecewise-constant empirical accelerations\cite{beutler_celestialTheoretical_2010,beutler_celestialApplication_2010}. 

In Kalman filtering, a common ad hoc method is to manually tune a diagonal process noise covariance. Such an approach does not scale the process noise covariance according to the length of the propagation interval, so the covariance must be retuned for different propagation interval lengths\cite{carpenter_navigation_2018}. Furthermore, this technique does not capture the cross-covariance between elements of the spacecraft state, and does not provide any dynamical consistency between the magnitudes of the process noise covariance of the different spacecraft state components. A second approach in Kalman filtering is to model the effect on the spacecraft state error covariance of specific forces. For example, the effects of auto-correlated central body gravity modeling errors of both commission and omission can be modeled\cite{schiemenz_improved_2020,wright_sequential_1981}. However, these models are rarely used onboard due to computational cost\cite{carpenter_navigation_2018}. The process noise due to maneuver execution errors can also be explicitly modeled\cite{carpenter_navigation_2018,sullivan_nonlinear_2020}. Olson et al.\cite{olson_precomputing_2017} describe how consider covariance analysis techniques can be used to precompute a process noise covariance profile along a reference trajectory to capture the effects of known model parameter uncertainties on the spacecraft state error covariance. Although this approach is highly accurate provided that the assumed parameter uncertainties are representative and the actual spacecraft trajectory is close to the predefined reference trajectory, the numerical precomputation process is computationally expensive. Another approach called dynamic model compensation directly estimates the unmodeled accelerations by augmenting the estimated state vector with empirical accelerations\cite{cruickshank_genetic_1998,schutz_statistical_2004,leonard_gravity_2013,myers_dynamical_1975,tapley_orbit_1973}. Augmenting the state vector increases filter computation time, and the performance of dynamic model compensation depends on the tuning of parameters in the dynamical model of the empirical accelerations such as the empirical acceleration time correlation constants\cite{goldstein_real-time_2001,leonard_gravity_2013}.

Advances in process noise modeling have primarily focused on an absolute Cartesian spacecraft state. However, the most appropriate state representation depends on the application. The spacecraft state can be parameterized as absolute or relative Cartesian coordinates as well as absolute or relative orbital elements. An absolute state describes the motion of a single spacecraft. In contrast, a relative state is a function of the absolute states of two spacecraft and describes the motion of one spacecraft relative to the other. Cartesian coordinates are most frequently used because of their more direct relation to the modeled measurements and because they do not have any singularities. Although orbital element states suffer from singularities, they have several advantages over Cartesian coordinates. Since most orbital elements vary more slowly in time than position and velocity, orbital elements can often be numerically integrated more efficiently\cite{sullivan_angles_2018}. Orbital elements also enable analytical perturbation models\cite{guffanti_linear_2019,koenig_new_2017} and provide greater geometric insight than a Cartesian state. See Hintz\cite{hintz_survey_2008} and Sullivan et al.\cite{sullivan_comprehensive_2017} for comprehensive reviews of absolute and relative state parameterizations respectively.

One of the most widely used process noise models in sequential orbit determination is state noise compensation (SNC)\cite{schutz_statistical_2004,carpenter_navigation_2018,myers_filtering_1974,shalom_estimation_2002}. Although the unmodeled accelerations are generally correlated in time to some degree, SNC treats the unmodeled accelerations as zero-mean Gaussian white process noise. The resulting process noise covariance can be evaluated through numerical integration, which is computationally intensive. Alternatively, an analytical SNC model is commonly used that assumes kinematic motion with zero nominal spacecraft acceleration. However, this model is restricted to an absolute Cartesian state and is only valid for small propagation intervals\cite{carpenter_navigation_2018,myers_filtering_1974,shalom_estimation_2002}. The authors are not aware of analytical SNC models for absolute orbital element state representations or mathematically rigorous, analytical process noise covariance models for relative spacecraft states.

To overcome these limitations, this paper develops new analytical SNC process noise covariance models for both absolute and relative spacecraft states. The following section reviews relevant background information on SNC and spacecraft relative motion dynamical modeling. Exploiting relative dynamics models, the subsequent section develops new analytical process noise covariance models for absolute Cartesian and orbital element states by assuming either a circular orbit or a small propagation interval. These models are then extended to long propagation intervals for perturbed, eccentric orbits. The next two sections present frameworks for modeling the process noise covariance of relative spacecraft states by assuming either small or large interspacecraft separation. The large interspacecraft separation framework is very general, and can leverage any absolute process noise covariance model, not just SNC models. The presented analytical models are validated through numerical simulations in the penultimate section. Finally, conclusions are provided based on the numerical results.

\section{Background}
\subsection{State Noise Compensation for an Absolute Spacecraft State}
Let $\bm{x}_{\alpha}(t) \in \mathbb{R}^{6}$ denote the absolute state of a spacecraft parameterized in either Cartesian coordinates or a set of orbital elements. The specific state representation is indicated by the subscript $\alpha$. The mean state estimate at time step $i$ after processing all the measurements through time step $j$ is $\hat{\bm{x}}_{\alpha,i|j}$. Assuming an unbiased estimator such that $\text{E}[\hat{\bm{x}}_{\alpha,i|j}] = \bm{x}_{\alpha}(t_i)$, the error covariance or formal covariance of the mean state estimate at time step $i$ after processing all the measurements through time step $j$ is $\mathbf{P}_{\alpha,i|j} = \text{E}[(\bm{x}_{\alpha}(t_i) - \hat{\bm{x}}_{\alpha,i|j})(\bm{x}_{\alpha}(t_i) -\hat{\bm{x}}_{\alpha,i|j})^T]$. An essential task in sequential orbit determination as well as satellite conjunction analysis is to propagate the error covariance to some future time $t_k$, which involves determining $\mathbf{P}_{\alpha,k|k-1}$ given $\hat{\bm{x}}_{\alpha,k-1|k-1}$ and $\mathbf{P}_{\alpha,k-1|k-1}$. Note that $t_k$ can be any time of interest greater than $t_{k-1}$.

In general, the dynamical model of the spacecraft state is given by some nonlinear function $\dot{\bm{x}}_{\alpha}(t) = \bm{f}(\bm{x}_{\alpha}(t),\bm{u}(t),t)$ where $\bm{u}(t)$ is the control input. The approach taken in a linearized framework such as an extended Kalman filter is to linearize the dynamical model over the propagation interval $[t_{k-1},t_k]$ about the mean state estimate taking into account all measurements through time $t_{k-1}$.  This linearization results in a linear time-varying system
\begin{equation}\label{eq:ct LTV dynamical model}
	\dot{\bm{X}}_{\alpha}(t) = \mathbf{A}_{\alpha}(t)\bm{X}_{\alpha}(t) + \mathbf{B}_{\alpha}(t)\bm{u}(t) + \mathbf{\Gamma}_{\alpha}(t)\bm{\epsilon}(t)
\end{equation}
where $\bm{X}_{\alpha}(t) = \bm{x}_{\alpha}(t) - \hat{\bm{x}}_{\alpha}(t)$ is the state estimate error, and $\hat{\bm{x}}_{\alpha}(t)$ is the mean state estimate at time $t$ taking into account all measurements through time $t_{k-1}$ such that $\hat{\bm{x}}_{\alpha}(t_k) = \hat{\bm{x}}_{\alpha,k|k-1}$. Here $\mathbf{A}_{\alpha}$ is the plant matrix, $\mathbf{B}_{\alpha}$ is the control input matrix, and $\mathbf{\Gamma}_{\alpha}$ is the process noise mapping matrix. These matrices are defined by the partial derivatives
\begin{equation}
\begin{aligned}
  \mathbf{A}_{\alpha}(t) &= \frac{\partial \dot{\bm{x}}_{\alpha}(t)}{\partial \bm{x}_{\alpha}(t)},   
  \hspace{.2cm}&   
  \mathbf{B}_{\alpha}(t) &= \frac{\partial \dot{\bm{x}_{\alpha}}(t)}{\partial \bm{u}(t)}, 
  \hspace{.2cm}& 
  \mathbf{\Gamma}_{\alpha}(t) &= \frac{\partial \dot{\bm{x}}_{\alpha}(t)}{\partial \bm{\epsilon}(t)}\\
\end{aligned}
\end{equation}
evaluated at the mean state estimate $\hat{\bm{x}}_{\alpha}(t)$.

The continuous-time process noise $\bm{\epsilon}$ describes stochastic deviations from the nominal spacecraft dynamics. If the dynamics were truly linear, the only sources of process noise would be numerical error and unmodeled accelerations. Throughout this paper $\bm{\epsilon}$ physically represents unmodeled accelerations since they generally create orders of magnitude more process noise than numerical error. Here $\bm{\epsilon} \in \mathbb{R}^{3}$ is modeled as a zero-mean white Gaussian process with autocovariance 
\begin{equation}\label{eq:absolute power spectral density}
	\text{E}[\bm{\epsilon}(t)\bm{\epsilon}(\tau)^T] 	= \mathbf{\widetilde{Q}}\delta(t-\tau)	
\end{equation}
where $\delta(\cdot)$ is the Dirac delta function and $\mathbf{\widetilde{Q}}\in \mathbb{R}^{3\times 3}$ is the process noise power spectral density, which describes the strength of the unmodeled accelerations. This approach to process noise modeling is called SNC\cite{schutz_statistical_2004,carpenter_navigation_2018}. The stochastic vector of unmodeled accelerations $\bm{\epsilon}$ can be inertial or radial-transverse-normal (RTN) coordinates, which will be denoted by the superscripts $\mathcal{I}$ and $\mathcal{R}$ respectively. Thus $\bm{\epsilon}^{\mathcal{I}}(t) = \ $\tiny$\raisebox{1.8pt}{$\underset{\mathcal{R}\rightarrow \mathcal{I}}{}$}$  \hspace{-.5cm}\small $\raisebox{2pt}{$\mathbf R$}$ \hspace{.01cm} \normalsize$(t) \bm{\epsilon}^\mathcal{R}(t)$ where \mbox{\tiny$\raisebox{1.8pt}{$\underset{\mathcal{R}\rightarrow \mathcal{I}}{}$}$  \hspace{-.5cm}\small $\raisebox{2pt}{$\mathbf R$}$ \hspace{.01cm} \normalsize(t) $\in \mathbb{R}^{3\times 3}$} is the rotation matrix from the RTN frame to the inertial frame at time $t$. The power spectral density of $\bm{\epsilon}^\mathcal{I}$ and $\bm{\epsilon}^\mathcal{R}$ are $\widetilde{\mathbf{Q}}^{\mathcal{I}}$ and $\widetilde{\mathbf{Q}}^{\mathcal{R}}$ respectively. In the RTN frame, also referred to as the Hill frame, the radial and normal axes are aligned with the chief spacecraft position and angular momentum vectors respectively. The transverse axis completes the right-handed triad and is positive in the direction of the chief velocity. Typically, the spacecraft dynamics are highly nonlinear, and $\mathbf{\widetilde{Q}}$ should be large enough to also accommodate errors due to dynamical nonlinearities in the propagation of the mean state estimate and associated error covariance.

The discrete-time solution of Eq. (\ref{eq:ct LTV dynamical model}) is
\begin{equation}\label{eq:dt LTV dynamical model}
	\bm{X}_{\alpha}(t_k) = \mathbf{\Phi}_{\alpha}(t_k,t_{k-1})\bm{X}_{\alpha}(t_{k-1}) + \int_{t_{k-1}}^{t_k} \mathbf{\Phi}_{\alpha}(t_k,\tau)\mathbf{B}_{\alpha}(\tau)\bm{u}(\tau) d\tau  + \bm{w}_{\alpha,k}
\end{equation}
where the Jacobian $\mathbf{\Phi}_{\alpha}(t_k,t_{k-1}) = \frac{\partial \bm{x}_{\alpha}(t_k)}{\partial \bm{x}_{\alpha}(t_{k-1})}$ evaluated at the mean state estimate is the state transition matrix. The discrete-time process noise is $\bm{w}_{\alpha,k} \sim \mathcal{N}(\bm{0},\mathbf{Q}_{\alpha,k})$, and the process noise covariance is
\begin{equation}\label{eq:absolute process noise covariance}
 	\mathbf{Q}_{\alpha,k} = \int_{t_{k-1}}^{t_k}  \mathbf{\Phi}_{\alpha}(t_k,\tau)  \mathbf{\Gamma}_{\alpha}(\tau)  \mathbf{\widetilde{Q}}  \mathbf{\Gamma}_{\alpha}(\tau)^T  \mathbf{\Phi}_{\alpha}(t_k,\tau)^T  d\tau
\end{equation}
Due to the structure of Eq. (\ref{eq:absolute process noise covariance}), $\mathbf{Q}_{\alpha,k}$ is guaranteed symmetric and positive semi-definite provided that $\mathbf{\widetilde{Q}}$ is symmetric and positive semi-definite. The propagated or time-updated error covariance is
\begin{align}
		\mathbf{P}_{\alpha,k|k-1} &= \text{E}[\bm{X}_{\alpha}(t_{k}) \bm{X}_{\alpha}(t_{k})^T]\\
									  &=\mathbf{\Phi}_{\alpha,k}\mathbf{P}_{\alpha,k-1|k-1}\mathbf{\Phi}_{\alpha,k}^T+\mathbf{Q}_{\alpha,k}\label{eq:covariance time update}
\end{align}
where $\mathbf{\Phi}_{\alpha,k} = \mathbf{\Phi}_{\alpha}(t_k,t_{k-1})$.

The process noise covariance integral in Eq. (\ref{eq:absolute process noise covariance}) can be evaluated numerically provided that the state transition matrix and process noise mapping matrix are integrable functions. However, analytical solutions are desirable because they significantly reduce computation time and provide insight into system behavior.
One common analytical approximation of Eq. (\ref{eq:absolute process noise covariance}) considers an absolute, inertial Cartesian state 
\begin{equation}\label{eq:inertial Cart coordinates}
	\bm{x}_{\mathcal{I}} = 
	\begin{bmatrix}
		\bm{r}\\ 
		\bm{v}
	\end{bmatrix}
\end{equation}
where $\bm{r}$ and $\bm{v}$ are the inertial position and velocity vectors respectively. Assuming kinematic motion where the nominal spacecraft acceleration is zero, 
\begin{align}\label{eq:STM kinematic}
	\mathbf{\Phi}_{\mathcal{I}}(t_{k},t_{k-1}) = 
	\begin{bmatrix}
		\mathbf{I}_{3\times 3} 	&\Delta t_{k}\mathbf{I}_{3\times 3}\\
		\mathbf{0}_{3\times 3}	&\mathbf{I}_{3\times 3}\\
	\end{bmatrix}
\end{align}
where $\Delta t_{k} = t_{k} - t_{k-1}$ is the length of the propagation interval $[t_{k-1},t_k]$. The identity matrix and matrix of zeros with three rows and columns are $\mathbf{I}_{3\times 3}$ and $\mathbf{0}_{3\times 3}$ respectively. Modeling $\bm{\epsilon}$ in the inertial frame, the process noise mapping matrix is 
\begin{equation}\label{eq:Gamma Cart}
	\mathbf{\Gamma}_{\mathcal{I}} = 
	\begin{bmatrix}
		\mathbf{0}_{3\times 3}\\
		\mathbf{I}_{3\times 3}
	\end{bmatrix}
\end{equation}
Substituting Eqs. (\ref{eq:STM kinematic}-\ref{eq:Gamma Cart}) and $\mathbf{\widetilde{Q}} = \mathbf{\widetilde{Q}}^{\mathcal{I}}$ into Eq. (\ref{eq:absolute process noise covariance}) and then evaluating the integral yields\cite{carpenter_navigation_2018,shalom_estimation_2002}
\begin{align}\label{eq:Q kinematic inertial}
    \mathbf{Q}_{\mathcal{I},k} = 
     \begin{bmatrix}
        \frac{1}{3}\Delta t_k^3 \mathbf{\widetilde{Q}}^{\mathcal{I}}  &\frac{1}{2}\Delta t_k^2 \mathbf{\widetilde{Q}}^{\mathcal{I}}\\
        \frac{1}{2}\Delta t_k^2\mathbf{\widetilde{Q}}^{\mathcal{I}}   &\Delta t_k\mathbf{\widetilde{Q}}^{\mathcal{I}}
    \end{bmatrix}
\end{align}
Assuming the orientation of the RTN frame relative to the inertial frame is constant over a small propagation interval $[t_{k-1},t_k]$, the relation \mbox{$\widetilde{\mathbf{Q}}^{\mathcal{I}} = \ $\tiny$\raisebox{1.8pt}{$\underset{\mathcal{R}\rightarrow \mathcal{I}}{}$}$  \hspace{-.5cm}\small $\raisebox{2pt}{$\mathbf R$}$ \hspace{.01cm} \normalsize$   \widetilde{\mathbf{Q}}^{\mathcal{R}}$\tiny$\raisebox{1.8pt}{$\underset{\mathcal{R}\rightarrow \mathcal{I}}{}$}$  \hspace{-.5cm}\small $\raisebox{2pt}{$\mathbf R$}$ \hspace{.01cm} \normalsize$^T$} can be substituted into Eq. (\ref{eq:Q kinematic inertial})\cite{carpenter_navigation_2018}. Another approach to simplify Eq. (\ref{eq:absolute process noise covariance}) is to model $\bm{\epsilon}$ as constant over the interval $[t_{k-1},t_k]$\cite{schutz_statistical_2004}. However, the covariance of $\bm{\epsilon}$ must then be retuned for different propagation interval lengths as observed by Carpenter et. al\cite{carpenter_navigation_2018}. 
Although the analytical process noise covariance model in Eq. (\ref{eq:Q kinematic inertial}) is widely used, it is only valid for small propagation intervals because it assumes kinematic motion with zero nominal spacecraft acceleration\cite{carpenter_navigation_2018,shalom_estimation_2002}. Furthermore, the authors are not aware of  any analytical SNC models for orbital element states or mathematically rigorous, analytical process noise covariance models for relative spacecraft states. These gaps in the state of the art are addressed in this paper.

\subsection{Hill-Clohessy-Wiltshire Equations}
\begin{figure}[!h]
\centering
\includegraphics[width=.85\linewidth,trim=205 45 280 165,clip]{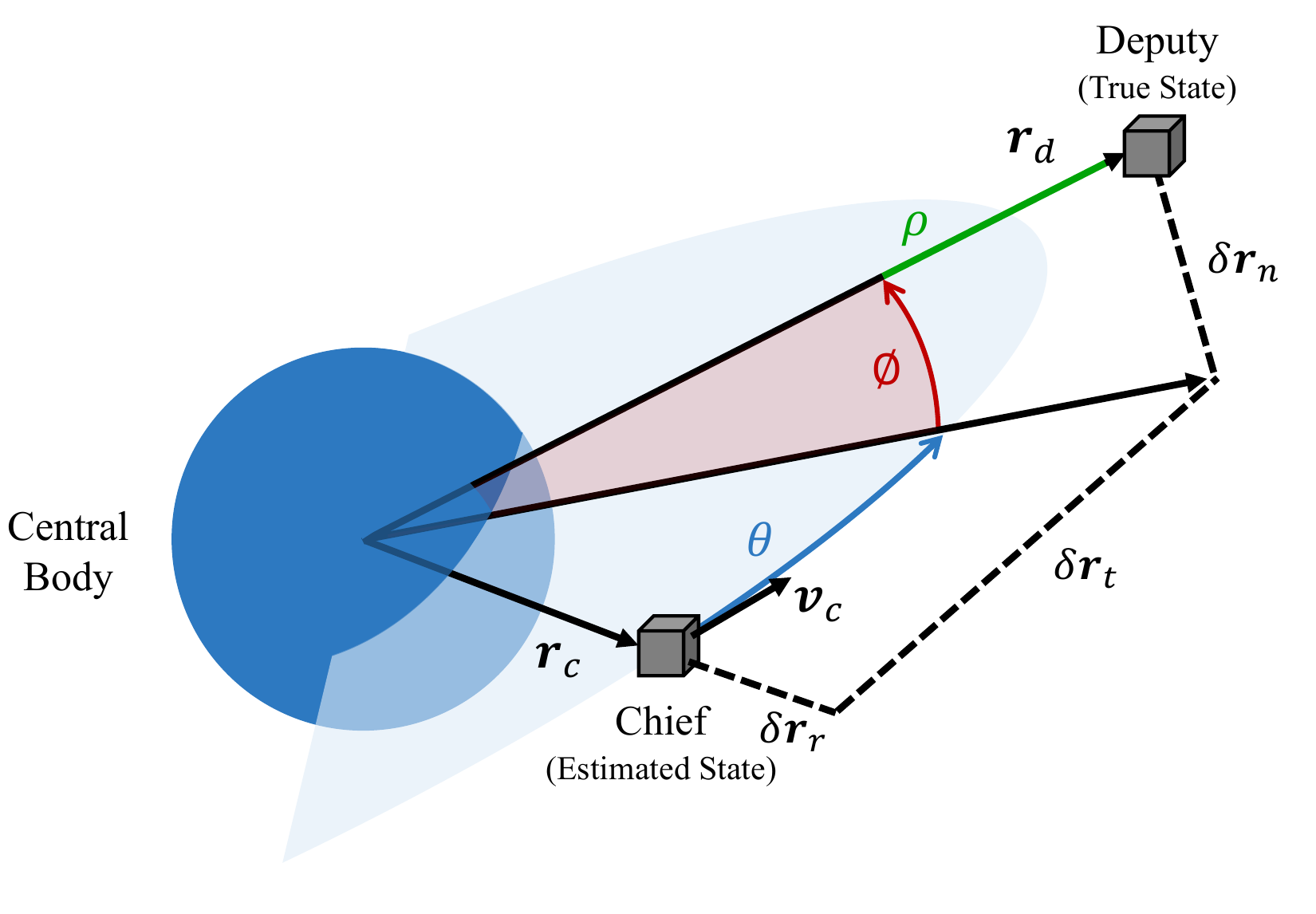}
\caption{Rectilinear and curvilinear relative spacecraft coordinates. It is illustrated that the chief and deputy can be considered the estimated and true absolute states respectively in order to model the error dynamics of an absolute spacecraft state.}
\label{fig:relative coordinates}
\end{figure}

Many spacecraft relative dynamics models have been developed to describe the motion of a deputy spacecraft relative to a chief spacecraft\cite{sullivan_comprehensive_2017,koenig_new_2017,guffanti_linear_2019}. Any relative spacecraft dynamics model can be used to model the dynamics of the state error by considering the estimated and true spacecraft states to be the chief and deputy spacecraft respectively as illustrated in Figure \ref{fig:relative coordinates}. One particularly simple relative dynamics model is the Hill-Clohessy-Wiltshire (HCW) equations, which assume two-body motion, a circular chief orbit, and small interspacecraft separation as compared to the chief orbit radius. The solution of the HCW equations was used by Geller\cite{geller_orbital_2007} to approximate a scalar metric of the uncertainty of the position of one spacecraft relative to another spacecraft due to navigation and maneuver execution errors as well as umodeled accelerations. In the following sections, the solution of the HCW equations is used for the first time to construct an analytical model of the full process noise covariance of both absolute and relative Cartesian states. This section reviews the HCW equations\cite{clohessy_terminal_1960}, which can be parameterized in either rectilinear or curvilinear coordinates, and the corresponding discrete-time solution\cite{alfriend_spacecraft_2010,de_comparative_2011}.

\vspace{.5cm}


Throughout this paper, $\delta \bm{x}_{\alpha}(t) \in \mathbb{R}^{6}$ denotes a relative spacecraft state where the specific relative state representation is denoted by the subscript $\alpha$. Let
\begin{equation}\label{eq:relative cartesian state}
\delta\bm{x}_{\mathcal{R}} =
	\begin{bmatrix}
		\delta \bm{r}\\
		\delta \bm{v}
	\end{bmatrix}
	=
	\begin{bmatrix}
		\underset{\mathcal{I} \rightarrow \mathcal{R}}{\mathbf R}(\bm{r}_d - \bm{r}_c)\\
		\underset{\mathcal{I} \rightarrow \mathcal{R}}{\mathbf R}(\bm{v}_d - \bm{v}_c) - \bm{w}^{\times}\underset{\mathcal{I} \rightarrow \mathcal{R}}{\mathbf R}(\bm{r}_d - \bm{r}_c)
	\end{bmatrix}
\end{equation}
be the state of a deputy spacecraft relative to a chief expressed in the chief RTN frame. The inertial Cartesian coordinates, defined in Eq. (\ref{eq:inertial Cart coordinates}), of the chief and deputy are denoted by the subscripts $c$ and $d$ respectively. The chief orbit may be the orbit of either a physical or virtual spacecraft. Notice that for the relative velocity vector $\delta \bm{v}$, the time derivative has been taken in the chief RTN frame. Here $\delta \bm{r} = [\delta r_r \ \ \delta r_t \ \ \delta r_n]^T$ and $\delta \bm{v} = [\delta v_r \ \ \delta v_t \ \ \delta v_n]^T$. The angular velocity vector of the chief RTN frame with respect to the inertial frame expressed in the chief RTN frame is denoted $\bm{w}$. The cross product matrix of $\bm{w}$ is $\bm{w}^{\times}$, which is a skew-symmetric matrix defined such that
\begin{equation}\label{eq:w cross product mat}
	\bm{w} = 
	\begin{bmatrix}
		w_r\\
		w_t\\
		w_n
	\end{bmatrix}
	\implies
	\bm{w}^{\times} = 
	\begin{bmatrix}
		0  		&-w_n 		&w_t\\
		w_n 		&0 			&-w_r\\
	    -w_t		&w_r 		&0
	\end{bmatrix}
\end{equation}
For two-body motion, $w_r = w_t = 0$ and $w_n = \dot{\theta_c}$ where $\dot{\theta_c}$ is the time rate of change of the angle between the chief position vector $\bm{r}_c$ and an arbitrary fixed vector in the chief orbital plane.

Consider the spherical coordinates
\begin{align}
	\rho &= ((r_c+\delta r_r)^2 + \delta r_t^2 + \delta r_n^2)^{1/2} - r_c\\
	\theta &= \text{atan}\left(\frac{\delta r_{t}}{r_{c} + \delta r_r}\right)\\
	\phi &= \text{atan}\left(\frac{\delta r_{n}}{((r_c + \delta r_{r})^2+\delta r_{t}^2)^{1/2}}  \right)
\end{align}
where $r_c = ||\bm{r}_c||$. Here $\rho = ||\bm{r}_d|| - r_c$, $\theta$ is the angle between the projection of the deputy position vector onto the chief orbital plane and the chief position vector, and $\phi$ is the angle between the deputy position vector and its projection onto the chief orbital plane. The angles $\theta$ and $\phi$ are measured positive in the directions of the chief velocity and angular momentum vectors respectively as shown in Figure \ref{fig:relative coordinates}. The nonlinear equations of relative motion parameterized using $\rho$, $\theta$, and $\phi$ are derived by Willis et al.\cite{willis_secondorder_2019} considering only two-body motion. Linearizing these equations about zero interspacecraft separation, assuming a circular chief orbit, and including differential perturbing accelerations yields\cite{willis_secondorder_2019,alfriend_spacecraft_2010}
\begin{align}
	\ddot{\rho} &=  3n^2 \rho + 2n r_c\dot{\theta} + \delta d_r \label{eq:curv eom linearized small t 1}\\
r_c\ddot{\theta} &= - 2n\dot{\rho} + \delta d_t\\
r_c\ddot{\phi} &= -n^2r_c\phi + \delta d_n \label{eq:curv eom linearized small t 3}
\end{align}
The chief mean motion is denoted $n$, and $\delta\bm{d} = \bm{d}_d - \bm{d}_c = [\delta d_r \ \ \delta d_t \ \ \delta d_n]^T$ is the difference between the perturbing accelerations of the deputy and chief expressed in the chief RTN frame. 

One particular curvilinear relative state is
\begin{equation}\label{eq:rel curv state}
	\delta\bm{x}_\psi = [\rho \ \ r_c \theta \ \ r_c \phi \ \ \dot{\rho} \ \ r_c \dot{\theta} \ \ r_c\dot{\phi}]^T
\end{equation}
The discrete-time solution of Eqs. (\ref{eq:curv eom linearized small t 1}-\ref{eq:curv eom linearized small t 3}) yields the state transition matrix 
\begin{equation}\label{eq:hcw stm}
	\mathbf{\Phi}_{\delta \psi}(t,0) = 
	\begin{bmatrix}
		4-3\cos nt 	&0	&0	&\frac{1}{n}\sin nt 	&\frac{2}{n}(1-\cos nt) 	&0\\
		6(\sin nt -nt) 	&1	&0	&\frac{2}{n}(\cos nt-1) 	&\frac{1}{n}(4\sin nt-3nt) 	&0\\
		0	&0 	&\cos nt 	&0 	&0 	&\frac{1}{n}\sin nt\\
		3n\sin nt 	&0 	&0 	&\cos nt 	&2\sin nt 	&0\\
		6n(\cos nt -1)	&0 	&0 	&-2\sin nt 	&4\cos nt-3 	&0\\
		0	&0	&-n\sin nt 	&0 	&0 	&\cos nt
	\end{bmatrix}
\end{equation}
that propagates $\delta\bm{x}_\psi$ from time $0$ to $t$, which is the well-known solution of the HCW equations. The HCW equations can also be parameterized in the rectilinear coordinates $\delta\bm{x}_\mathcal{R}$ defined in Eq. (\ref{eq:relative cartesian state}), and the corresponding state transition matrix is identical to that shown in Eq. (\ref{eq:hcw stm}). Parameterizing the state in curvilinear coordinates instead of rectilinear coordinates has been shown to reduce errors incurred by the linearization when using the state transition matrix to propagate the relative state\cite{alfriend_spacecraft_2010,de_comparative_2011}.

\section{Absolute Process Noise Covariance}
In general, it is difficult to obtain an exact analytical solution of the integral in Eq. (\ref{eq:absolute process noise covariance}) that is reasonably concise even when only considering two-body motion. For an absolute Cartesian state, the complexity of the integral in Eq. (\ref{eq:absolute process noise covariance}) is principally due to the state transition matrix. The predominant approach in the literature is to simplify the Cartesian state transition matrix by assuming kinematic motion with zero nominal acceleration, which results in the analytical model in Eq. (\ref{eq:Q kinematic inertial}) that is valid for small propagation intervals. An alternate approach is developed in the following subsection by assuming a circular orbit and leveraging the solution of the HCW equations. The resulting analytical process noise covariance model is accurate for much longer propagation intervals for near-circular orbits than the widely used model in Eq. (\ref{eq:Q kinematic inertial}). In contrast to a Cartesian state, the complexity of the integral in Eq. (\ref{eq:absolute process noise covariance}) for an orbital element state is primarily due to the process noise mapping matrix, which is given by the Gauss Variational Equations. New analytical process noise covariance models for absolute orbital element states are derived in the following subsections by simplifying the process noise mapping matrix through the assumption of either a circular orbit or a short propagation interval. These absolute Cartesian and orbital element process noise covariance models are then extended to long propagation intervals for perturbed, eccentric orbits by exploiting spacecraft relative dynamics models.



Typically, it is advantageous to model spacecraft acceleration uncertainty in the RTN frame since the amount of acceleration uncertainty in each axis of the RTN frame generally varies less than in each axis of an inertial frame. For example, consider a case where atmospheric drag is the largest contributor to spacecraft acceleration uncertainty. Throughout the orbit, the acceleration uncertainty is consistently large in the along-track direction and relatively small in the radial and normal directions. In contrast, the level of acceleration uncertainty in each axis of an inertial frame will change significantly throughout the orbit. Consequently, the unmodeled acceleration power spectral density, $\widetilde{\mathbf{Q}}$, will be modeled in the RTN frame throughout the rest of this paper. For simplicity, it is assumed that the power spectral density matrix is diagonal such that 
\begin{equation}\label{eq:psd RTN}
	\widetilde{\mathbf{Q}}^\mathcal{R} = 
	\begin{bmatrix}
		\widetilde{Q}_r 	&0 							&0\\
		0 							&\widetilde{Q}_t		&0\\
		0 							&0 							&\widetilde{Q}_n
	\end{bmatrix}
\end{equation}
Here $\widetilde{Q}_r$, $\widetilde{Q}_t$, and $\widetilde{Q}_n$ describe the strength of the unmodeled accelerations in the radial, transverse, and normal directions respectively. A diagonal $\widetilde{\mathbf{Q}}^\mathcal{R}$ enables the expression for $\mathbf{Q}_{\alpha,k}$ in Eq. (\ref{eq:absolute process noise covariance}) to be written as
\begin{equation}\label{eq:Q split RTN}
	\mathbf{Q}_{\alpha,k} = \mathbf{Q}^{r}_{\alpha,k} + \mathbf{Q}^{t}_{\alpha,k} + \mathbf{Q}^{n}_{\alpha,k}
\end{equation}
Here
\begin{equation}\label{eq:Q in each axis}
	\mathbf{Q}^i_{\alpha,k} = \widetilde{Q}_i \int^{t_k}_{t_{k-1}} \overline{\mathbf{\Gamma}}_{\alpha,i}(t_k,\tau) \overline{\mathbf{\Gamma}}_{\alpha,i}(t_k,\tau)^T d\tau
\end{equation}
is the contribution to $\mathbf{Q}_{\alpha,k}$ due to unmodeled accelerations in axis $i$ of the RTN frame where $i\in \{r,t,n \}$. The vectors $\overline{\mathbf{\Gamma}}_{\alpha,r}(t_k,\tau)$, $\overline{\mathbf{\Gamma}}_{\alpha,t}(t_k,\tau)$, and $\overline{\mathbf{\Gamma}}_{\alpha,n}(t_k,\tau)$ are the first, second, and third columns respectively of
\begin{equation}\label{eq:Gamma bar}
	\overline{\mathbf{\Gamma}}_{\alpha}(t_k,\tau) = \mathbf{\Phi}_{\alpha}(t_k,\tau)  \mathbf{\Gamma}_{\alpha}(\tau)
\end{equation}
The particular values of $\widetilde{Q}_r$, $\widetilde{Q}_t$, and $\widetilde{Q}_n$ are based on a priori knowledge of the dynamical environment. For example, these values can be selected in order to match elements of $\mathbf{Q}_{\alpha,k}$ to corresponding estimates\cite{shalom_estimation_2002} or in order for the filter error covariance to match an empirical approximation of the true error covariance\cite{carpenter_navigation_2018}. Alternatively, the power spectral density can be estimated online using the recently developed adaptive state noise compensation algorithm\cite{stacey_adaptive_2021}.

\subsection{Cartesian State and Circular Orbit}
This section leverages the solution of the HCW equations to develop a process noise covariance model for the absolute, inertial Cartesian state of a spacecraft $\bm{x}_{\mathcal{I}}$ for near-circular orbits considering only two-body motion. The state error of an inertial Cartesian state $\bm{X}_{\mathcal{I}}(t) = \bm{x}_{\mathcal{I}}(t) - \hat{\bm{x}}_{\mathcal{I}}(t)$ can be thought of as the state of a deputy relative to a chief spacecraft. The deputy is the true spacecraft state $\bm{x}_{\mathcal{I}}$, and the chief is the estimated state $\hat{\bm{x}}_{\mathcal{I}}$, which is a virtual orbit (See Figure \ref{fig:relative coordinates}). This relative motion is described in a linearized framework through Eqs. (\ref{eq:ct LTV dynamical model}) and (\ref{eq:dt LTV dynamical model}). The relative motion can be equivalently parameterized using the curvilinear state $\delta \bm{x}_{\psi}$ defined in Eq. (\ref{eq:rel curv state}). The parameterization of $\bm{X}_{\mathcal{I}}$ in curvilinear coordinates will be denoted $\delta \bm{X}_{\psi}$. Assuming an unbiased estimator, the expectation of the estimate error is zero (i.e., $\text{E}[\delta \bm{X}_{\psi}(t)] = \bm{0}$). The dynamics of $\delta\bm{X}_\psi$ can be linearized about $\delta \bm{X}_\psi = \bm{0}$, resulting in a linear time-varying system similar to Eq. (\ref{eq:ct LTV dynamical model}). Following the same approach used to derive Eqs. (\ref{eq:dt LTV dynamical model}-\ref{eq:covariance time update}), the covariance of $\delta\bm{X}_\psi$ is
\begin{align}
		\mathbf{P}_{\psi,k|k-1} &= \text{E}[\delta\bm{X}_{\psi}(t_{k}) \delta\bm{X}_{\psi}(t_{k})^T]\\
									  &=\mathbf{\Phi}_{\psi,k}\mathbf{P}_{\psi,k-1|k-1}\mathbf{\Phi}_{\psi,k}^T+\mathbf{Q}_{\psi,k}\label{eq:covariance time update curv}
\end{align}
Here $\mathbf{\Phi}_{\psi,k}$ is the state transition matrix of $\delta\bm{X}_{\psi}$. Assuming two-body motion and a circular orbit, $\mathbf{\Phi}_{\psi,k}$ is equivalent to Eq. (\ref{eq:hcw stm}). In this case, $n$ in Eq. (\ref{eq:hcw stm}) refers to the mean motion computed from the mean state estimate. The unmodeled accelerations $\bm{\epsilon}$ are modeled in the RTN frame of the mean state estimate trajectory. The process noise mapping matrix is then deduced from Eqs. (\ref{eq:curv eom linearized small t 1}-\ref{eq:curv eom linearized small t 3}) to be
\begin{equation}\label{eq:Gamma Cart Rot}
	\mathbf{\Gamma}_{\mathcal{\psi}} = 
	\begin{bmatrix}
		\mathbf{0}_{3\times 3}\\
		\mathbf{I}_{3\times 3}
	\end{bmatrix}
\end{equation}

The matrix $\mathbf{Q}_{\psi,k}$ is obtained by substituting Eqs. (\ref{eq:hcw stm}-\ref{eq:psd RTN}) and Eq. (\ref{eq:Gamma Cart Rot}) into Eq. (\ref{eq:absolute process noise covariance}).
The solution to Eq. (\ref{eq:absolute process noise covariance}) is computed by evaluating Eq. (\ref{eq:Q in each axis}) and substituting the result into Eq. (\ref{eq:Q split RTN}). Analytically evaluating Eq. (\ref{eq:Q in each axis}) for each axis yields
\begin{flalign}\label{eq:Q HCW components rel rad}
\mathbf{Q}^{r}_{\psi,k} &= \widetilde{Q}_r
		\begin{bmatrix}
			\frac{1}{n^2}\zeta_{ss} 					&\ddots															&\ddots 	&\ddots 		&\ddots 		&\ddots\\
			\frac{2}{n^2}(\zeta_{cs} - \zeta_s) 	&\frac{4}{n^2}(\zeta_{cc}-2\zeta_c+\Delta t_k) 		&\ddots 	&\ddots 		&\ddots 		&\ddots\\
			0 													&0 																	&0 			&\ddots 		&\ddots 		&\ddots\\
			\frac{1}{n}\zeta_{cs} 						&\frac{2}{n}(\zeta_{cc} - \zeta_c) 						&0 			&\zeta_{cc} 	&\ddots 		&\ddots\\
			 \frac{-2}{n}\zeta_{ss}					&\frac{-4}{n}(\zeta_{cs}-\zeta_s)			 				&0 			&-2\zeta_{cs} &4\zeta_{ss} 	&\ddots\\
			0 													&0 																	&0 			&0 				&0 				&0
		\end{bmatrix}&
\end{flalign}
\begin{flalign}\label{eq:Q HCW components rel tra}
\mathbf{Q}^{t}_{\psi,k} &= \widetilde{Q}_t
\left[\begin{matrix}
\frac{4}{n^2}(\Delta t_k+\zeta_{cc}-2\zeta_c)  							&\ddots 	&\ddots \\
			\frac{2}{n^2}(3n\zeta_{tc}-4\zeta_{cs}+4\zeta_s-\frac{3n}{2}\Delta t_k^2)  	&\frac{1}{n^2}(3n^2\Delta t_k^3+16\zeta_{ss}-24n\zeta_{ts}) 	&\ddots\\
			0  	&0 	&0\\
			\frac{4}{n}(\zeta_s-\zeta_{cs})  	&\frac{2}{n}(4\zeta_{ss}-3n\zeta_{ts}) 	&0\\
			\frac{2}{n}(7\zeta_c-4\zeta_{cc}-3\Delta t_k)  	&\frac{1}{n}(16\zeta_{cs}-12n\zeta_{tc}-12\zeta_s+\frac{9n}{2}\Delta t_k^2) 	&0\\
			0  	&0 	&0
\end{matrix}\nonumber \right.\\
&\qquad\qquad \qquad \qquad \qquad \qquad \qquad \qquad \quad \ \ \
\left.\begin{matrix}
			\ddots 	&\ddots 	&\ddots\\
			\ddots 	&\ddots 	&\ddots\\
			\ddots 	&\ddots 	&\ddots\\
			4\zeta_{ss} 	&\ddots 	&\ddots\\
			8\zeta_{cs}-6\zeta_s 	&9\Delta t_k+16\zeta_{cc}-24\zeta_c 	&\ddots\\
			0 	&0 	&0
\end{matrix}\right]&
\end{flalign}
\begin{flalign}\label{eq:Q HCW components rel nor}
\mathbf{Q}^{n}_{\psi,k} &= \widetilde{Q}_n
		\begin{bmatrix}
			0  	&\ddots 	&\ddots 						&\ddots 	&\ddots 	&\ddots\\
			0  	&0 			&\ddots 						&\ddots 	&\ddots 	&\ddots\\
			0  	&0 			&\frac{1}{n^2}\zeta_{ss}&\ddots 	&\ddots 	&\ddots\\
			0  	&0 			&0 								&0 			&\ddots 	&\ddots\\
			0  	&0			&0 								&0 			&0 			&\ddots\\
			0  	&0 			&\frac{1}{n}\zeta_{cs} 	&0 			&0 			&\zeta_{cc}\\
		\end{bmatrix}&
\end{flalign}
Each matrix in Eqs. (\ref{eq:Q HCW components rel rad}-\ref{eq:Q HCW components rel nor}) is symmetric, and only the unique elements in the lower triangular portion of each matrix are specified for brevity. The auxiliary variables in Eqs. (\ref{eq:Q HCW components rel rad}-\ref{eq:Q HCW components rel nor}) are
\begingroup
\allowdisplaybreaks
\begin{equation}\label{eq:HCW integrals relative}
\begin{aligned}
	\zeta_c &= \int_{t_{k-1}}^{t_k} \text{cos}(n(t_k-\tau)) \ d\tau = \frac{1}{n}\sin n\Delta t_k\\
	\zeta_s &= \int_{t_{k-1}}^{t_k} \text{sin}(n(t_k-\tau)) \ d\tau = \frac{1}{n}(1- \cos n\Delta t_k)\\
	\zeta_{cc} &= \int_{t_{k-1}}^{t_k} \text{cos}^2(n(t_k-\tau)) \ d\tau = \frac{\Delta t_k}{2} + \frac{1}{4n}\sin  2n\Delta t_k\\
	\zeta_{ss} &= \int_{t_{k-1}}^{t_k} \text{sin}^2(n(t_k-\tau)) \ d\tau = \frac{\Delta t_k}{2} - \frac{1}{4n}\sin  2n\Delta t_k\\
	\zeta_{cs} &= \int_{t_{k-1}}^{t_k} \text{cos}(n(t_k-\tau))\text{sin}(n(t_k-\tau)) \ d\tau = \frac{1}{4n}(1 - \cos 2n\Delta t_k)\\
	\zeta_{tc} &= \int_{t_{k-1}}^{t_k} (t_k-\tau)\text{cos}(n(t_k-\tau)) \ d\tau = \frac{1}{n^2}(\cos n\Delta t_k - 1) + \frac{\Delta t_k}{n} \sin n\Delta t_k\\
		\zeta_{ts} &= \int_{t_{k-1}}^{t_k} (t_k-\tau)\text{sin}(n(t_k-\tau)) \ d\tau = \frac{-\Delta t_k}{n}\cos n\Delta t_k + \frac{1}{n^2}\sin n\Delta t_k
\end{aligned}
\end{equation}
\endgroup


Recalling the assumption of an unbiased filter such that $\text{E}[\bm{X}_{\mathcal{I}}(t_k)] = \text{E}[\delta\bm{X}_{\psi}(t_k)] = \bm{0}$, $\bm{X}_{\mathcal{I}}(t_k)$ can be written in terms of its curvilinear parameterization through a first-order Taylor series expansion about $\bm{X}_{\mathcal{I}}(t_k) = \bm{0}$ as
\begin{equation}\label{eq:taylor series X_I}
	\bm{X}_{\mathcal{I}}(t_k) = \mathbf{J}_{\mathcal{I}\psi}(t_k) \delta\bm{X}_\psi(t_k)
\end{equation}
where
\begin{equation}\label{eq:partials error wrt rectilinear}
	\mathbf{J}_{\mathcal{I}\psi}(t_k) 
	= \frac{\partial \delta\bm{x}_{\mathcal{I}}(t_k)}{\partial\delta \bm{x}_\psi(t_k)}\bigg\rvert_{\delta\bm{x}_\psi = \bm{0}}  	
\end{equation}
Here $\delta\bm{x}_{\mathcal{I}} = \bm{x}_{\mathcal{I}}^d - \bm{x}_{\mathcal{I}}^c$ is the difference between the absolute, inertial Cartesian states of a deputy and chief. Recall that in this case the chief and deputy are the estimated and true states respectively (See Figure \ref{fig:relative coordinates}). Eq. (\ref{eq:partials error wrt rectilinear}) can be expanded using the chain rule as
\begin{equation}\label{eq:JIpsi chain}
	\mathbf{J}_{\mathcal{I}\psi}(t_k) = \left(\frac{\partial \delta \bm{x}_{\mathcal{I}}(t_k)}{\partial \ \delta\bm{x}_{\mathcal{R}}(t_k)}\bigg\rvert_{\delta\bm{x}_{\mathcal{R}}(t_k)=\bm{0}}\right)        \left(\frac{\partial \ \delta\bm{x}_{\mathcal{R}}(t_k)}{\partial \ \delta\bm{x}_\psi(t_k)}\bigg\rvert_{\delta \bm{x}_\psi(t_k) = \bm{0}}\right) 
\end{equation}
 The first matrix of partial derivatives in Eq. (\ref{eq:JIpsi chain}) is deduced from Eq. (\ref{eq:relative cartesian state}) to be 
\begin{align}\label{eq:partials dxI wrt dxR}
	\frac{\partial \delta\bm{x}_{\mathcal{I}}(t_k)}{\partial \ \delta\bm{x}_{\mathcal{R}}(t_k)} \bigg\rvert_{\delta\bm{x}_{\mathcal{R}}(t_k)=\bm{0}}
	= 
	\begin{bmatrix}
		\underset{\mathcal{R}\rightarrow \mathcal{I}}{\mathbf R}(t_k) 		&\mathbf{0}_{3 \times 3}\\
		\underset{\mathcal{R}\rightarrow \mathcal{I}}{\mathbf R}(t_k)\bm{w}^{\times}   		&\underset{\mathcal{R}\rightarrow \mathcal{I}}{\mathbf R}(t_k)
	\end{bmatrix}
\end{align}
Since the relative Cartesian and curvilinear coordinates are equivalent to first order at zero separation for a circular orbit,
$\frac{\partial \ \delta\bm{x}_{\mathcal{R}}(t_k)}{\partial \ \delta\bm{x}_\psi(t_k)}\big\rvert_{\delta\bm{x}_\psi(t_k) = \bm{0}} = \mathbf{I}_{6\times 6}$. 

Using the linearization in Eq. (\ref{eq:taylor series X_I}), the time updated error covariance of the inertial Cartesian state can be written as
\begingroup
\begin{align}
	\mathbf{P}_{\mathcal{I},k|k-1} &= \text{E}[\bm{X}_{\mathcal{I}}(t_k)\bm{X}_{\mathcal{I}}^T(t_k)]\label{eq:P f(QR) abs 1}\\
		&= \mathbf{J}_{\mathcal{I}\psi}(t_k) \mathbf{P}_{\psi,k|k-1} \mathbf{J}_{\mathcal{I}\psi}(t_k) ^T\label{eq:P f(QR) abs 1b}\\
								  &= \mathbf{J}_{\mathcal{I}\psi}(t_k)  \mathbf{\Phi}_{\psi,k} \mathbf{P}_{\psi,k-1|k-1} \mathbf{\Phi}^T_{\psi,k} \mathbf{J}_{\mathcal{I}\psi}(t_k) ^T + \mathbf{J}_{\mathcal{I}\psi}(t_k)  \mathbf{Q}_{\psi,k} \mathbf{J}_{\mathcal{I}\psi}(t_k)^T\label{eq:P f(QR) abs 2}\\				
								  &= \mathbf{\Phi}_{\mathcal{I},k}\mathbf{P}_{\mathcal{I},k-1|k-1}  \mathbf{\Phi}^T_{\mathcal{I},k} + \mathbf{J}_{\mathcal{I}\psi}(t_k)  \mathbf{Q}_{\psi,k} \mathbf{J}_{\mathcal{I}\psi}(t_k) ^T\label{eq:P f(QR) abs 3}
\end{align}
\endgroup
Substituting Eq. (\ref{eq:covariance time update curv}) into Eq. (\ref{eq:P f(QR) abs 1b}) yields Eq. (\ref{eq:P f(QR) abs 2}). Eqs. (\ref{eq:P f(QR) abs 2}) and (\ref{eq:P f(QR) abs 3}) are related through the chain rule and the linearization in Eq. (\ref{eq:taylor series X_I}) applied at time $t_{k-1}$. Comparing Eq. (\ref{eq:P f(QR) abs 3}) and Eq. (\ref{eq:covariance time update}), it is clear that the process noise covariance of $\bm{x}_{\mathcal{I}}$ at time step $k$ is 
\begin{equation}\label{eq:QI from Qpsi}
	\mathbf{Q}_{\mathcal{I},k} = \mathbf{J}_{\mathcal{I}\psi}(t_k) \mathbf{Q}_{\psi,k} \mathbf{J}_{\mathcal{I}\psi}(t_k)^T
\end{equation}
Since $\frac{\partial \ \delta\bm{x}_{\mathcal{R}}(t_k)}{\partial \ \delta\bm{x}_\psi(t_k)}\big\rvert_{\delta\bm{x}_\psi(t_k) = \bm{0}} = \mathbf{I}_{6\times 6}$, it is equivalent to compute the process noise covariance integral parameterized in the rectlinear coordinates $\delta\bm{x}_\mathcal{R}$, and then obtain $\mathbf{Q}_{\mathcal{I},k}$ through a linear mapping similar to Eq. (\ref{eq:QI from Qpsi}). Thus a higher order mapping from $\mathbf{Q}_{\psi,k}$ to $\mathbf{Q}_{\mathcal{I},k}$ is required to glean the benefits of spherical coordinates over rectilinear coordinates.


\subsection{Orbital Element State and Circular Orbit}
This section describes a framework for modeling the process noise covariance of an absolute orbital element state by assuming a circular orbit. The set of equinoctial elements is considered as an example, but the developed approach can be applied to any orbital element state representation. The equinoctial elements are defined as
\begin{equation} \label{eq:ee}
	\bm{x}_{\mathcal{E}} = 
	\begin{bmatrix}
	a\\
	f\\
	g\\
	h\\
	k\\
	\lambda
	\end{bmatrix}
	=
	\begin{bmatrix}
	a\\
	e\text{cos}(\omega+\Omega)\\
	e\text{sin}(\omega+\Omega)\\
	\text{tan}(\frac{i}{2})\text{cos}(\Omega)\\
	\text{tan}(\frac{i}{2})\text{sin}(\Omega)\\
	M+\omega+\Omega
	\end{bmatrix}
\end{equation}
in terms of the classical Keplerian orbital elements $a$, $e$, $i$, $\Omega$, $w$, and $M$. The state $\bm{x}_{\mathcal{E}}$ is singular for $i = 180 \degree$ and $e=1$. Considering only two-body motion, the state transition matrix of $\bm{x}_{\mathcal{E}}$ is
\begin{equation}\label{eq:stm oe}
\mathbf{\Phi}_{\mathcal{E}}(t,t_0) = 
\begin{bmatrix}
	1 								&\bm{0}_{1\times 4}		&\bm{0}\\
	\bm{0}_{4\times 1}	&\mathbf{I}_{4\times 4}	&\bm{0}_{4\times 1}\\
	\frac{-3n}{2a}(t-t_{0})&\bm{0}_{1\times 4} 	&1
\end{bmatrix}
\end{equation}
The time derivatives of the equinoctial elements are given by the Gauss Variational Equations, which can be written in matrix form as
\begin{equation}\label{eq:abslute GVE}
	\frac{d \bm{x}_{\mathcal{E}}}{dt} = \mathbf{G}(\bm{x}_{\mathcal{E}})\bm{d} + 
	\begin{bmatrix}
	\bm{0}_{5\times 1} \\
	n
	\end{bmatrix}
\end{equation}
where $\bm{d} \in \mathbb{R}^{3}$ is the perturbing accelerations expressed in the RTN frame. Treating the process noise as unmodeled accelerations expressed in the RTN frame, $\mathbf{G}(\bm{x}_{\mathcal{E}}) \in \mathbb{R}^{6\times 3}$ is the process noise mapping matrix. This matrix is\cite{battin_introduction_1987,roth_gaussian_1985}
\begin{align}\label{eq:Gamma equinoctial}
	\mathbf{\Gamma}_{\mathcal{E}} = \mathbf{G}(\bm{x}_{\mathcal{E}})
			=
			\begin{bmatrix}
				\bar{A}		&\bar{B} 		&0\\
				\bar{C}		&\bar{D}		&\bar{E}\\
				\bar{F}		&\bar{G}		&\bar{H}\\
				0				&0				&\bar{I}\\
				0 				&0				&\bar{J}\\
				\bar{K}		&\bar{L}		&\bar{M}
			\end{bmatrix}
\end{align}
where $\bar{A},\bar{B},\dots,\bar{M}$ are each a scalar function of $\bm{x}_{\mathcal{E}}$, and their specific definitions are provided in Appendix A. The matrix $\mathbf{\Gamma}_{\mathcal{E}}$ contains trigonometric functions of the true longitude $l = \nu + \Omega + w$ where $\nu$ is the true anomaly. For a near-circular orbit, the approximations $f=g=0$ and 
\begin{equation}\label{eq:true longitude circular}
	l(t) = l(t_{k-1}) + \bar{n}_k(t-t_{k-1})
\end{equation} 
can be made where 
\begin{equation}\label{eq:average angular rate}
	\bar{n}_k = \frac{1}{\Delta t_k} \int_{t_{k-1}}^{t_k} \dot{\theta}_s(t) \ dt = \frac{\Delta \theta_{s,k}}{\Delta t_k}
\end{equation}
is the average spacecraft angular rate over the propagation interval $[t_{k-1},t_k]$. These approximations become exact as the orbit eccentricity approaches zero. Here $\dot{\theta}_s$ is the time rate of change of the angle between the spacecraft position vector and an arbitrary fixed vector in the orbital plane. The variable $\Delta \theta_{s,k} = \theta_s(t_k) - \theta_s(t_{k-1})$ is the angle traversed by the spacecraft over the interval $[t_{k-1},t_k]$. In the case of zero eccentricity, $\bar{n}_k$ is equal to the mean motion. For other orbital element state representations, any quickly varying orbital elements such as the true anomaly or argument of latitude can be approximated similar to Eq. (\ref{eq:true longitude circular}).

The process noise covariance of the state $\bm{x}_{\mathcal{E}}$ at time step $k$, $\mathbf{Q}_{\mathcal{E},k}$, is computed by evaluating Eq. (\ref{eq:Q in each axis}) and substituting into Eq. (\ref{eq:Q split RTN}). Analytically evaluating Eq. (\ref{eq:Q in each axis}) for each axis using Eq. (\ref{eq:true longitude circular}) and assuming $f=g=0$ yields
\begingroup
\allowdisplaybreaks
\begin{align}
\mathbf{Q}^r_{\mathcal{E},k} &=
	\widetilde{Q}_r\left(\frac{p}{\mu}\right)
	\begin{bmatrix}
		0 		&\ddots 		&\ddots		&\ddots	&\ddots	&\ddots\\
		0 		&\bar{\zeta}_{ss}	&\ddots		&\ddots	&\ddots	&\ddots\\
		0 		&-\bar{\zeta}_{cs} &\bar{\zeta}_{cc}	&\ddots	&\ddots	&\ddots\\
		0 		&0 				&0				&0			&\ddots	&\ddots\\
		0 		&0 				&0				&0			&0			&\ddots\\
		0 		&-2\bar{\zeta}_s		&2\bar{\zeta}_{c}	&0			&0			&4\Delta t_k\\
	\end{bmatrix}\label{eq:QR equinoctial}\\
\mathbf{Q}^t_{\mathcal{E},k} &=
	\widetilde{Q}_t(\mathbf{M}_t\mathbf{M}_t^T) \circ
	\begin{bmatrix}
		\Delta t_k 					&\ddots 		&\ddots		&\ddots	&\ddots	&\ddots\\
		\bar{\zeta}_{c} 						&\bar{\zeta}_{cc}	&\ddots		&\ddots	&\ddots	&\ddots\\
		\bar{\zeta}_s							&\bar{\zeta}_{cs} 	&\bar{\zeta}_{ss}	&\ddots	&\ddots	&\ddots\\
		0 									&0 				&0				&0			&\ddots	&\ddots\\
		0 									&0 				&0				&0			&0			&\ddots\\
		\frac{1}{2}\Delta t_k^2 &\bar{\zeta}_{tc}	&\bar{\zeta}_{ts}	&0			&0			&\frac{1}{3}\Delta t_k^3\\
	\end{bmatrix}\\
\mathbf{Q}^n_{\mathcal{E},k} &=
	\widetilde{Q}_n(\mathbf{M}_n\mathbf{M}_n^T) \circ
	\begin{bmatrix}
		0 		&\ddots 	&\ddots	&\ddots							&\ddots							&\ddots\\
		0 		&0			&\ddots	&\ddots							&\ddots							&\ddots\\
		0 		&0 			&0			&\ddots							&\ddots							&\ddots\\
		0 		&0 			&0			&\bar{\zeta}_{cc}						&\ddots							&\ddots\\
		0 		&0 			&0			&\bar{\zeta}_{cs}						&\bar{\zeta}_{ss}						&\ddots\\
		0 		&0			&0			&k\bar{\zeta}_{cc}-h\bar{\zeta}_{cs}	&k\bar{\zeta}_{cs}-h\bar{\zeta}_{ss}	&k^2\bar{\zeta}_{cc}+h^2\bar{\zeta}_{ss}-2kh\bar{\zeta}_{cs}\\
	\end{bmatrix}\label{eq:QN equinoctial}
\end{align}
\endgroup
Each matrix in Eqs. (\ref{eq:QR equinoctial}-\ref{eq:QN equinoctial}) is symmetric, so only the unique elements in the lower triangular portion of each matrix are specified. The symbol $\circ$ denotes the Hadamard product, which indicates element-wise multiplication of two matrices with the same dimensions. The auxilliary variables in Eqs. (\ref{eq:QR equinoctial}-\ref{eq:QN equinoctial}) are
\begin{align}
\mathbf{M}_t &= \left[\frac{2a^2}{L} \ \ \ 2\sqrt{\frac{p}{\mu}} \ \ \ 2\sqrt{\frac{p}{\mu}} \ \ \ 0 \ \ \ 0 \ \ \ \frac{-3na}{L}\right]^T\\
\mathbf{M}_n &= \left[0 \ \ \ 0 \ \ \ 0 \ \ \ \frac{1}{2}\sqrt{\frac{p}{\mu}}(1+h^2+k^2) \ \ \ \frac{1}{2}\sqrt{\frac{p}{\mu}}(1+h^2+k^2) \ \ \ \frac{-L}{\mu}\right]^T
\end{align}
and
\begingroup
\allowdisplaybreaks
\begin{equation}\label{eq:Q equinoctial integration}
\begin{aligned}
	\bar{\zeta}_c 	&= \int_{t_{k-1}}^{t_k} \cos l(\tau)d\tau  =  \frac{1}{\bar{n}_k}(\sin l_{k} - \sin l_{k-1})\\
	\bar{\zeta}_s 	&= \int_{t_{k-1}}^{t_k} \sin l(\tau)d\tau  =  -\frac{1}{\bar{n}_k}(\cos l_{k} - \cos l_{k-1})\\
	\bar{\zeta}_{cc} &= \int_{t_{k-1}}^{t_k} \text{cos}^2 \ l(\tau)d\tau  = \frac{\Delta t_k}{2} + \frac{1}{4\bar{n}_k}(\sin 2l_{k} - \sin 2l_{k-1})\\
	\bar{\zeta}_{ss} &= \int_{t_{k-1}}^{t_k} \text{sin}^2 \ l(\tau)d\tau  = \frac{\Delta t_k}{2} - \frac{1}{4\bar{n}_k}(\sin 2l_{k} - \sin 2l_{k-1})\\
	\bar{\zeta}_{cs} &= \int_{t_{k-1}}^{t_k} \cos l(\tau) \ \sin l(\tau)d\tau = -\frac{1}{4\bar{n}_k}(\cos  2l_{k} - \cos  2l_{k-1})\\
	\bar{\zeta}_{tc} &= \int_{t_{k-1}}^{t_k} (t_k-\tau)\cos l(\tau) d\tau= -\frac{1}{\bar{n}_k^2}(\bar{n}_k\Delta t_k\sin l_{k-1} + \cos l_{k} - \cos l_{k-1})\\
	\bar{\zeta}_{ts} &= \int_{t_{k-1}}^{t_k} (t_k-\tau)\sin l(\tau)d\tau = \frac{1}{\bar{n}_k^2}(\bar{n}_k\Delta t_k\cos l_{k-1} - \sin l_{k} + \sin l_{k-1})
\end{aligned}
\end{equation}
\endgroup
Here $p = a(1-e^2)$ is the orbit semi-parameter, $L = \sqrt{\mu p}$ is the magnitude of the specific angular momentum vector, and $l_{k} = l(t_{k})$. Since the orbit is assumed circular, the semi-parameter can be replaced with the semi-major axis. 

\subsection{Orbital Element State and Small Propagation Interval}
A concise analytical approximation of Eq. (\ref{eq:absolute process noise covariance}) that is valid for eccentric orbits can be obtained by assuming that the propagation interval $[t_{k-1},t_k]$ is small such that any quickly varying orbital elements can be approximated as constant in $\mathbf{\Gamma}_{\mathcal{\alpha}}$ over that interval. Then $\mathbf{\Gamma}_{\mathcal{\alpha}}$ is constant when considering only two-body motion. In the case of the equinoctial elements, it is assumed that the true longitude is constant over the interval $[t_{k-1},t_k]$. For a constant $\mathbf{\Gamma}_{\mathcal{E}}$ and the state transition matrix in Eq. (\ref{eq:stm oe}), the contributions to $\mathbf{Q}_{\mathcal{E},k}$ in each axis shown in Eq. (\ref{eq:Q in each axis}) evaluate to
\begin{align}
	\mathbf{Q}_{\mathcal{E},k}^{r} &= \widetilde{Q}_r\left(\Delta t_k\mathbf{\Gamma}_{\mathcal{E},r}\mathbf{\Gamma}_{\mathcal{E},r}^{T} + \frac{3n\Delta t_k^2}{4a}
	\begin{bmatrix}
		 \bm{0}_{5\times5} 	&\mathbf{S}_r\\
		 \mathbf{S}_r^T  		&\frac{n}{a} \bar{A}^2\Delta t - 2\bar{A}\bar{K}
	\end{bmatrix}
	\right)\nonumber
	\\
	\mathbf{Q}_{\mathcal{E},k}^{t} &= \widetilde{Q}_t\left(\Delta t_k\mathbf{\Gamma}_{\mathcal{E},t}\mathbf{\Gamma}_{\mathcal{E},t}^T + \frac{3n\Delta t_k^2}{4a}
	\begin{bmatrix}
		 \bm{0}_{5\times5} 	&\mathbf{S}_{t}\\
		\mathbf{S}_t^T 			&\frac{n}{a} \bar{B}^2\Delta t - 2\bar{B}\bar{L}
	\end{bmatrix} 
	\right)\label{eq:Q equi small interval}
	\\
	\mathbf{Q}_{\mathcal{E},k}^{n} &= \widetilde{Q}_n\Delta t_k\mathbf{\Gamma}_{\mathcal{E},n}\mathbf{\Gamma}_{\mathcal{E},n}^T \nonumber
\end{align}
where $\mathbf{S}_r = - [\bar{A}^2 \ \ \bar{A}\bar{C} \ \ \bar{A}\bar{F} \ \ 0 \ \ 0]^T$ and $\mathbf{S}_t = -[\bar{B}^2 \ \ \bar{B}\bar{D} \ \ \bar{B}\bar{G} \ \ 0 \ \ 0]^T$. The first, second, and third columns of $\mathbf{\Gamma}_{\mathcal{E}}$ are denoted $\mathbf{\Gamma}_{\mathcal{E},r}$, $\mathbf{\Gamma}_{\mathcal{E},t}$, and $\mathbf{\Gamma}_{\mathcal{E},n}$ respectively. Eq. (\ref{eq:Q equi small interval}) is evaluated using the true longitude halfway through the propagation interval $[t_{k-1},t_k]$ at time $(t_{k-1}+t_k)/2$. The true longitude at time $(t_{k-1}+t_k)/2$ can be approximated through two-body propagation of the true longitude from time $t_{k-1}$.




Both of the developed equinoctial element process noise covariance models provide important insights. For example, notice from Eqs. (\ref{eq:abslute GVE}-\ref{eq:Gamma equinoctial}) that the time evolution of the equinoctial elements $h$ and $k$ are only affected by perturbing accelerations in the normal direction. As a result, only unmodeled accelerations in the normal direction contribute to the process noise covariance of $h$ and $k$ as is shown in the model described by Eqs. (\ref{eq:QR equinoctial}-\ref{eq:Q equinoctial integration}) and the model in Eq. (\ref{eq:Q equi small interval}). In both of these models, the process noise covariance of $h$ is small when the true longitude is near $l = \frac{(2n-1)\pi}{2}$ throughout the interval $[t_{k-1},t_k]$ for any integer $n$. This occurs because perturbing accelerations have no effect on $h$ at $l = \frac{(2n-1)\pi}{2}$ since $\bar{I}$ in Eq. (\ref{eq:Gamma equinoctial}) is zero. Similarly, the process noise covariance of $k$ is small when the true longitude is near $l = n\pi$ throughout the interval $[t_{k-1},t_k]$.



\subsection{Long Propagation Intervals for Perturbed, Eccentric Orbits}

In the case of a long propagation interval for an orbit that is eccentric or highly perturbed, the analytical models in the previous sections may not provide sufficient accuracy. This section describes a methodology for such cases. The propagation interval $[t_{k-1},t_k]$ is broken into $N$ subintervals $[\bar{t}_0,\bar{t}_1]$, $[\bar{t}_1,\bar{t}_2]$, $\dots$, $[\bar{t}_{N-1},\bar{t}_{N}]$ where $t_{k-1} = \bar{t}_0 < \bar{t}_1 < \dots < \bar{t}_{N} = t_{k}$. Using Eq. (\ref{eq:covariance time update}) to sequentially propagate the error covariance over each subinterval, it is evident that Eq. (\ref{eq:absolute process noise covariance}) can be equivalently written as
\begin{equation}\label{eq:Q pert eccen}
	\mathbf{Q}_{\alpha,k} = \sum_{i=1}^N \mathbf{\Phi}_{\alpha}(t_k,\bar{t}_i)  \mathbf{Q}_{\alpha}(\bar{t}_i,\bar{t}_{i-1})  \mathbf{\Phi}_{\alpha}(t_k,\bar{t}_i)^T
\end{equation}
where
\begin{equation}
 	\mathbf{Q}_{\alpha}(\bar{t}_i,\bar{t}_{i-1}) = \int_{\bar{t}_{i-1}}^{\bar{t}_i}  \mathbf{\Phi}_{\alpha}(\bar{t}_i,\tau)  \mathbf{\Gamma}_{\alpha}(\tau)  \mathbf{\widetilde{Q}}  \mathbf{\Gamma}_{\alpha}(\tau)^T  \mathbf{\Phi}_{\alpha}(\bar{t}_i,\tau)^T  d\tau
\end{equation}
is the process noise covariance computed over the subinterval $[\bar{t}_{i-1},\bar{t}_i]$. 

Formulating the process noise covariance according to Eq. (\ref{eq:Q pert eccen}) is advantageous because an analytical approximation that is valid over small propagation intervals can be used for each $\mathbf{Q}_{\alpha}(\bar{t}_i,\bar{t}_{i-1})$ provided that the length of each subinterval is sufficiently small. Then an analytical state transition matrix $\mathbf{\Phi}_{\alpha}(t_k,\bar{t}_i)$ with the accuracy required for the specified application can be used to map the effects of each $\mathbf{Q}_{\alpha}(\bar{t}_i,\bar{t}_{i-1})$ on $\mathbf{Q}_{\alpha,k}$. In this way, each term of the sum in Eq. (\ref{eq:Q pert eccen}) may neglect certain effects such as perturbations or orbit curvature over the subinterval $[\bar{t}_{i-1},\bar{t}_i]$, but these effects can be accounted for throughout the remainder of the propagation interval $[\bar{t}_i, t_k]$ through the selected $\mathbf{\Phi}_{\alpha}(t_k,\bar{t}_i)$. Furthermore, Eq. (\ref{eq:Q pert eccen}) provides the flexibility to use different values of $\mathbf{\widetilde{Q}}$ for different subintervals. Eq. (\ref{eq:Q pert eccen}) can also be used to extend many process noise covariance models from literature to perturbed, eccentric orbits such as the dynamic model compensation process noise covariance model developed by Cruickshank\cite{cruickshank_genetic_1998}, which assumes a small propagation interval\cite{carpenter_navigation_2018}.

Since $\mathbf{\Phi}_{\alpha}(t_k,\bar{t}_i)$ describes the evolution of the state estimate error (i.e., the true state relative to the mean state estimate), it can be formulated using any relative spacecraft state transition matrix (See Figure \ref{fig:relative coordinates}). 
Specifically, 
\begin{equation}\label{eq:stm from relative stm}
	\mathbf{\Phi}_{\alpha}(t_k,t_{k-1}) = \frac{\partial \delta \bm{x}_{\alpha}(t_k)}{\partial \delta \bm{x}_{\beta}(t_k)}  \mathbf{\Phi}_{\delta\beta}(t_k,t_{k-1})  \frac{\partial \delta \bm{x}_{\beta}(t_{k-1})}{\partial \delta \bm{x}_{\alpha}(t_{k-1})}
\end{equation}
where the partial derivatives are evaluated at zero interspacecraft separation, $\delta \bm{x}_{\alpha} = \bm{x}_{\alpha}^d - \bm{x}_{\alpha}^c$ is the difference between the states of a deputy and chief, and $\mathbf{\Phi}_{\delta\beta}(t_k,t_{k-1}) = \frac{\partial \delta \bm{x}_{\beta}(t_{k})}{\partial \delta \bm{x}_{\beta}(t_{k-1})}$ evaluated at zero separation is the state transition matrix of some relative state $\delta\bm{x}_{\beta}$. A comprehensive review of state transition matrices for relative states available in literature is provided by Sullivan et al.\cite{sullivan_comprehensive_2017} including models for orbits of any eccentricity and models that incorporate a variety or perturbations. Additional state transition matrices for relative orbital elements are provided by Koenig et al.\cite{koenig_new_2017} as well as Guffanti and D'Amico\cite{guffanti_linear_2019} that include the dominant effects of earth oblateness, atmospheric drag, and solar radiation pressure. In order to compute each $\mathbf{\Phi}_{\alpha}(t_k,\bar{t}_i)$ and $ \mathbf{Q}_{\alpha}(\bar{t}_i,\bar{t}_{i-1})$ in Eq. (\ref{eq:Q pert eccen}), the spacecraft state at each $\bar{t}_i$ is required. If the reference trajectory or mean state estimate is numerically integrated, the $\bar{t}_i$ can be chosen as times where the state is already available from the reference trajectory numerical integration.

As an example application of Eq. (\ref{eq:Q pert eccen}), consider a spacecraft in a near-circular orbit and a propagation interval that is multiple orbits. To compute each $\mathbf{Q}_{\alpha}(\bar{t}_i,\bar{t}_{i-1})$, Eqs. (\ref{eq:Q HCW components rel rad}-\ref{eq:HCW integrals relative}) and (\ref{eq:QI from Qpsi}) can be used for a Cartesian state or a model similar to Eqs. (\ref{eq:QR equinoctial}-\ref{eq:Q equinoctial integration}) can be used for an orbital element state. The subintervals $[\bar{t}_{i-1},\bar{t}_i]$ could be about one orbit or less in length since, as will be shown later, the effects of perturbations on the process noise covariance are relatively small over one orbit. Then an analytical $\mathbf{\Phi}_{\alpha}(t_k,\bar{t}_i)$ can be employed in Eq. (\ref{eq:Q pert eccen}) that accounts for perturbations to a sufficient degree to satisfy mission requirements.

As a second example, consider a spacecraft in an eccentric orbit and a propagation interval that is about one orbit. For an equinoctial element state, each $\mathbf{Q}_{\alpha}(\bar{t}_i,\bar{t}_{i-1})$ can be modeled through Eq. (\ref{eq:Q equi small interval}). Since the effects of perturbations are typically small over one orbit, the two-body state transition matrix in Eq. (\ref{eq:stm oe}) can be used for each $\mathbf{\Phi}_{\alpha}(t_k,\bar{t}_i)$. For an inertial Cartesian state, each $\mathbf{Q}_{\alpha}(\bar{t}_i,\bar{t}_{i-1})$ can be modeled through Eq. (\ref{eq:Q kinematic inertial}). The state transition matrix can be computed through Eq. (\ref{eq:stm from relative stm}) using any of the solutions of the Tschauner-Hempel equations such as the Yamanaka-Ankerson\cite{yamanaka_new_2002} state transition matrix. This matrix is formulated for a state $\delta \bm{x}_{\widetilde{\mathcal{R}}} = [\delta \widetilde{\bm{r}}^T \ \ \delta \widetilde{\bm{v}}^T]^T$ that is a normalized form of the relative state $\delta \bm{x}_{\mathcal{R}}$ defined in Eq. (\ref{eq:relative cartesian state}) where the time derivative is taken with respect to true anomaly. The transformations between these two states are\cite{willis_secondorder_2019}
\begin{equation}\label{eq:transformation RTN and normalized}
\begin{aligned}
\delta \widetilde{\bm{r}} &= \frac{1}{r_c} \delta \bm{r},   
  \hspace{.2cm}&   
 \delta \widetilde{\bm{v}} &= -\frac{e}{p} \delta \bm{r} \sin f + \frac{1}{1+e \cos f}\sqrt{\frac{p}{\mu}}\delta\bm{v}\\
  \delta \bm{r} &= r_c\delta \widetilde{\bm{r}},   
  \hspace{.2cm}&   
 \delta \bm{v} &= \sqrt{\frac{\mu}{p}}(e\delta\widetilde{\bm{r}} \sin f + (1+e \cos f )\delta\widetilde{\bm{v}})
\end{aligned}
\end{equation}
The required partial derivatives in Eq. (\ref{eq:stm from relative stm}) are
\begin{align}
	\frac{\partial \delta \bm{x}_{\mathcal{I}}(t_k)}{\partial \delta \bm{x}_{\widetilde{\mathcal{R}}}(t_k)} &= 
	\frac{\partial \delta \bm{x}_{\mathcal{I}}(t_k)}{\partial \delta \bm{x}_{\mathcal{R}}(t_k)} 
	\frac{\partial \delta \bm{x}_{\mathcal{R}}(t_k)}{\partial \delta \bm{x}_{\widetilde{\mathcal{R}}}(t_k)}\\
	\frac{\partial \delta \bm{x}_{\widetilde{\mathcal{R}}}(t_{k-1})}{\partial \delta \bm{x}_{\mathcal{I}}(t_{k-1})} &= 
	\frac{\partial \delta \bm{x}_{\widetilde{\mathcal{R}}}(t_{k-1})}{\partial \delta \bm{x}_{\mathcal{R}}(t_{k-1})}
	\frac{\partial \delta \bm{x}_{\mathcal{R}}(t_{k-1})}{\partial \delta \bm{x}_{\mathcal{I}}(t_{k-1})}
\end{align}
evaluated at zero separation where
\begin{align}
	\frac{\partial \delta \bm{x}_{\mathcal{R}}}{\partial \delta \bm{x}_{\widetilde{\mathcal{R}}}} &= 
	\begin{bmatrix}
		r_c\mathbf{I}_{3\times 3} 												    &\mathbf{0}_{3\times 3}\\
		e\sqrt{\frac{\mu}{p}}\sin f \ \mathbf{I}_{3\times 3}	&\sqrt{\frac{\mu}{p}}(1+e  \cos f ) \mathbf{I}_{3\times 3}
	\end{bmatrix}\label{eq:YA partials 1}\\
	\frac{\partial \delta \bm{x}_{\widetilde{\mathcal{R}}}}{\partial \delta \bm{x}_{\mathcal{R}}} &= 
	\begin{bmatrix}
		\frac{1}{r_c}\mathbf{I}_{3\times 3} 						&\mathbf{0}_{3\times 3}\\
		-\frac{e}{p}\sin f	\ \mathbf{I}_{3\times 3} 	&\frac{1}{1+e \cos f}\sqrt{\frac{p}{\mu}}\mathbf{I}_{3\times 3} 
	\end{bmatrix}\label{eq:YA partials 2}\\
	\frac{\partial \delta \bm{x}_{\mathcal{R}}}{\partial \delta \bm{x}_{\mathcal{I}}} &=
	\begin{bmatrix}
		\underset{\mathcal{I} \rightarrow \mathcal{R}}{\mathbf R} 								&\mathbf{0}_{3\times 3}\\
		- \bm{w}^{\times}\underset{\mathcal{I} \rightarrow \mathcal{R}}{\mathbf R}	&\underset{\mathcal{I} \rightarrow \mathcal{R}}{\mathbf R} 
	\end{bmatrix}\label{eq:partial dxR wrt dxI}
\end{align}
and $\frac{\partial \delta \bm{x}_{\mathcal{I}}}{\partial \delta \bm{x}_{\mathcal{R}}}$ is defined in Eq. (\ref{eq:partials dxI wrt dxR}). For near-circular orbits, true anomaly and argument of periapsis are poorly defined and can undergo abrupt changes due to perturbations. Thus, it is advantageous to use the Yamanaka-Ankerson state transition matrix formulated using the argument of latitude $u = w+f$ developed by Willis and D'Amico\cite{willis_analytical_2021}. The substitutions $e \cos  f = e_x\cos u + e_y\sin u$ and $e \sin f = e_x\sin u - e_y\cos u$ can also be made in Eqs. (\ref{eq:transformation RTN and normalized}) and (\ref{eq:YA partials 1}-\ref{eq:YA partials 2}) where the eccentricity vector components $e_x = e \cos w$ and $e_y = e \sin w$ are well defined for near-circular orbits. For both Cartesian coordinates and equinoctial elements, the accuracy of the resulting $\mathbf{Q}_{\alpha,k}$ increases as the length of the subintervals decreases since Eq. (\ref{eq:Q kinematic inertial}) and Eq. (\ref{eq:Q equi small interval}) assume a small propagation interval. Thus, the number of subintervals $N$ is chosen to balance accuracy and computational cost.

\section{Relative Process Noise Covariance For Small Separations}

If the interspacecraft separation is small compared to the orbit radius, then the nonlinear dynamical model of the relative state $\delta \bm{x}_{\alpha}$ can be linearized about zero interspacecraft separation over the propagation interval $[t_{k-1},t_k]$, resulting in a linear time-varying system
\begin{equation}\label{eq:curv ct LTV dynamical model}
	\delta\dot{\bm{x}}_{\alpha}(t) = \mathbf{A}_{\delta \alpha}(t)\delta\bm{x}_{\alpha}(t) + \mathbf{B}_{\delta \alpha}(t)\delta \bm{u}(t) + \mathbf{\Gamma}_{\delta\alpha}(t)\delta \bm{\epsilon}(t) 
\end{equation}
The plant matrix, control input matrix, and process noise mapping matrix of the relative state are
\begin{equation}
\begin{aligned}
  \mathbf{A}_{\delta\alpha}(t) &= \frac{\partial \delta\dot{\bm{x}}_{\alpha}(t)}{\partial \delta\bm{x}_{\alpha}(t)},   \hspace{.2cm}&   \mathbf{B}_{\delta\alpha}(t) &= \frac{\partial \delta\dot{\bm{x}}_{\alpha}(t)}{\partial \delta\bm{u}(t)}, \hspace{.2cm}& \mathbf{\Gamma}_{\delta\alpha}(t) &= \frac{\partial \delta\dot{\bm{x}}_{\alpha}(t)}{\partial \delta\bm{\epsilon}(t)}\\
\end{aligned}
\end{equation}
evaluated at zero interspacecraft separation. The subscript $\delta$ indicates a matrix corresponds to a relative state, and the particular relative state is specified by the subscript $\alpha$. Because the system is linearized about zero separation, the chief and deputy RTN frames are identical for the nominal or reference trajectory. For consistency, the RTN frame will be referred to as that of the chief throughout this section. The differential control inputs $\delta \bm{u} =  \bm{u}_d - \bm{u}_c$ are the difference between the control inputs of the deputy and chief expressed in the chief RTN frame. Similarly, the differential unmodeled accelerations $\delta \bm{\epsilon} =  \bm{\epsilon}_d - \bm{\epsilon}_c$ are the difference between the deputy and chief unmodeled accelerations expressed in the chief RTN frame with autocovariance 
\begin{align}
	\text{E}[\delta\bm{\epsilon}(t)\delta\bm{\epsilon}(\tau)^T] 
	&= \text{E}[\bm{\epsilon}_d(t) \bm{\epsilon}_d(\tau)^T] + \text{E}[\bm{\epsilon}_c(t) \bm{\epsilon}_c(\tau)^T] - \text{E}[\bm{\epsilon}_d(t) \bm{\epsilon}_c(\tau)^T] - \text{E}[\bm{\epsilon}_c(t) \bm{\epsilon}_d(\tau)^T]\nonumber\\
	&= (\mathbf{\widetilde{Q}}_d + \mathbf{\widetilde{Q}}_c - \mathbf{\widetilde{Q}}_{dc} - \mathbf{\widetilde{Q}}_{cd})\delta(t-\tau)\label{eq:cross power spectral density}\\
	&= \mathbf{\widetilde{Q}}_{\delta}\delta(t-\tau)\nonumber
\end{align}
Here $\mathbf{\widetilde{Q}}_{\delta} = (\mathbf{\widetilde{Q}}_d + \mathbf{\widetilde{Q}}_c - \mathbf{\widetilde{Q}}_{dc} - \mathbf{\widetilde{Q}}_{cd}) \in \mathbb{R}^{3\times 3}$ is the power spectral density of $\delta \bm{\epsilon}$. The power spectral densities of the chief and deputy unmodeled accelerations are $\mathbf{\widetilde{Q}}_c$ and $\mathbf{\widetilde{Q}}_d$ respectively. The cross power spectral densities are $\mathbf{\widetilde{Q}}_{dc}$ and $\mathbf{\widetilde{Q}}_{cd}$. Following a similar derivation to that of Eqs. (\ref{eq:absolute process noise covariance} - \ref{eq:covariance time update}), the propagated error covariance of $\delta \bm{x}_{\alpha}$ for the linearized system in Eq. (\ref{eq:curv ct LTV dynamical model}) is
\begin{align}
\mathbf{P}_{\delta\alpha,k|k-1} &= \text{E}[(\delta\bm{x}_{\alpha}(t_k) - \delta\hat{\bm{x}}_{\alpha,k|k-1})(\delta\bm{x}_{\alpha}(t_k) - \delta\hat{\bm{x}}_{\alpha,k|k-1})^T]\\
									  &=\mathbf{\Phi}_{\delta\alpha,k}\mathbf{P}_{\delta\alpha,k-1|k-1}\mathbf{\Phi}_{\delta\alpha,k}^T+\mathbf{Q}_{\delta\alpha,k}\label{eq:covariance time update relative}
\end{align}
where the process noise covariance of $\delta\bm{x}_{\alpha}$ at time step $k$ is
\begin{equation}\label{eq:relative process noise covariance}
 	\mathbf{Q}_{\delta\alpha,k} = \int_{t_{k-1}}^{t_k}  \mathbf{\Phi}_{\delta \alpha}(t_k,\tau)  \mathbf{\Gamma}_{\delta\alpha}(\tau)  \mathbf{\widetilde{Q}}_{\delta}  \mathbf{\Gamma}_{\delta\alpha}(\tau)^T  \mathbf{\Phi}_{\delta\alpha}(t_k,\tau)^T  d\tau
\end{equation}
The state transition matrix $\mathbf{\Phi}_{\delta\alpha,k} = \mathbf{\Phi}_{\delta\alpha}(t_k,t_{k-1})$ propagates the relative state $\delta \bm{x}_{\alpha}$ from time $t_{k-1}$ to time $t_k$. Due to the structure of Eq. (\ref{eq:relative process noise covariance}), $\mathbf{Q}_{\delta\alpha,k}$ is guaranteed symmetric and positive semi-definite when $\mathbf{\widetilde{Q}}_{\delta}$ is symmetric and positive semi-definite. Throughout this paper, $\mathbf{\widetilde{Q}}_{\delta}$ is expressed in the chief RTN frame, denoted $\widetilde{\mathbf{Q}}_{\delta}^\mathcal{R}$. For simplicity, it is assumed that $\widetilde{\mathbf{Q}}_{\delta}^\mathcal{R}$ is diagonal where
\begin{equation}\label{eq:Qtilde delta R}
	\widetilde{\mathbf{Q}}_{\delta}^\mathcal{R} = 
	\begin{bmatrix}
		\widetilde{Q}_{\delta r} 	&0 										&0\\
		0 										&\widetilde{Q}_{\delta t}	&0\\
		0 										&0 										&\widetilde{Q}_{\delta n}
	\end{bmatrix}
\end{equation}

\subsection{Relative Cartesian State}
This section models the process noise covariance of the relative state $\delta\bm{x}_{\mathcal{R}}$ defined in Eq. (\ref{eq:relative cartesian state}). First the process noise covariance is constructed for the relative state $\delta\bm{x}_{\mathcal{I}} = \bm{x}_{\mathcal{I}}^d - \bm{x}_{\mathcal{I}}^c$, which is the difference between the absolute, inertial Cartesian states of the deputy and chief. Linearizing the system about zero interspacecraft separation, the state transition matrix and process noise mapping matrix match those of the absolute inertial Cartesian state estimate error $\bm{X}_{\mathcal{I}}(t) = \bm{x}_{\mathcal{I}}(t) - \hat{\bm{x}}_{\mathcal{I}}(t)$. As a result, the expression for the process noise covariance of $\delta\bm{x}_{\mathcal{I}}$, denoted $\mathbf{Q}_{\delta \mathcal{I},k}$, is equivalent to $\mathbf{Q}_{\mathcal{I},k}$ except that $\mathbf{Q}_{\delta \mathcal{I},k}$ considers relative unmodeled accelerations. Thus, $\mathbf{Q}_{\delta \mathcal{I},k}$ can be approximated using Eq. (\ref{eq:Q kinematic inertial}) for a small propagation interval, Eqs. (\ref{eq:Q split RTN}) and (\ref{eq:Q HCW components rel rad}-\ref{eq:HCW integrals relative}) for a near-circular orbit, or Eq. (\ref{eq:Q pert eccen}) by making the substitutions $\mathbf{Q}_{\mathcal{I},k} = \mathbf{Q}_{\delta\mathcal{I},k}$ and $\widetilde{\mathbf{Q}}^\mathcal{R}= \widetilde{\mathbf{Q}}_{\delta}^\mathcal{R}$ where $\widetilde{\mathbf{Q}}_{\delta}^\mathcal{R}$ is defined in Eq. (\ref{eq:Qtilde delta R}). Then, similar to Eqs. (\ref{eq:P f(QR) abs 1}-\ref{eq:QI from Qpsi}), the process noise covariance of $\delta \bm{x}_{\mathcal{R}}$ is
\begin{equation}\label{eq:Qdr small separation}
	\mathbf{Q}_{\delta\mathcal{R},k}  =  \mathbf{J}_{\delta\mathcal{R}}(t_k) \mathbf{Q}_{\delta\mathcal{I},k} \mathbf{J}_{\delta\mathcal{R}}(t_k)^T
\end{equation}
where $\mathbf{J}_{\delta\mathcal{R}}(t_k) = \frac{\partial\delta \bm{x}_{\mathcal{R}}(t_k)}{\partial \delta\bm{x}_{\mathcal{I}}(t_k)}\bigg\rvert_{\delta\bm{x}_{\mathcal{I}}=\bm{0}}$ is given in Eq. (\ref{eq:partial dxR wrt dxI}).

\subsection{Relative Orbital Element State}
A framework for modeling the process noise covariance of a relative orbital element state is demonstrated using the set of relative orbital elements\cite{koenig_new_2017}
\begin{equation}\label{eq:roe equinoctial}
	\delta \bm{x}_{\mathcal{E}} = 
	\begin{bmatrix}
		(a_d - a_c)/a_c\\
		\lambda_d - \lambda_c\\ 
		f_d - f_c\\
		g_d - g_c\\
		h_d - h_c\\
		k_d - k_c
	\end{bmatrix}
\end{equation}
as an example. However, this approach can be applied to any set of relative orbital elements. The relative state in Eq. (\ref{eq:roe equinoctial}) is a function of the equinoctial orbital elements of the chief and deputy, which are defined in Eq. (\ref{eq:ee}). Consider a similar set of relative orbital elements
\begin{equation}
	\delta \bm{x}_{\mathcal{E}'} =  \bm{x}_{\mathcal{E}}^d - \bm{x}_{\mathcal{E}}^c
\end{equation}
where $\bm{x}_{\mathcal{E}}^c$ and $\bm{x}_{\mathcal{E}}^d$ are the equinoctial elements of the chief and deputy respectively. Linearizing the system about zero interspacecraft separation and considering only two-body motion, the state transition matrix of $\delta \bm{x}_{\mathcal{E}'}$ matches that shown in Eq. (\ref{eq:stm oe}). Since the system is linearized about zero interspacecraft separation, the semi-major axis and mean motion in the state transition matrix can refer to the chief or deputy. For consistency, the chief parameters will be used throughout this section whenever there is an arbitrary choice between the parameters of the chief and deputy.

The time derivative of $\delta \bm{x}_{\mathcal{E}'}$ can be expanded using the chain rule as
\begin{align}\label{eq:roe chain rule}
	\frac{d\delta \bm{x}_{\mathcal{E}'}}{dt} = \frac{\partial \delta \bm{x}_{\mathcal{E}'}}{\partial \bm{x}_{\mathcal{E}}^d} \frac{d \bm{x}_{\mathcal{E}}^d}{dt}  +  \frac{\partial \delta \bm{x}_{\mathcal{E}'}}{\partial \bm{x}_{\mathcal{E}}^c} \frac{d \bm{x}_{\mathcal{E}}^c}{dt} 
\end{align}
The partial derivatives of $\delta \bm{x}_{\mathcal{E}'}$ with respect to $\bm{x}_{\mathcal{E}}^d$ and $\bm{x}_{\mathcal{E}}^c$ are
\begin{equation}\label{eq:partials dxe dxed}
	\frac{\partial \delta \bm{x}_{\mathcal{E}'}}{\partial \bm{x}_{\mathcal{E}}^d} = \mathbf{I}_{6\times 6}, 
	\hspace{1cm} 
	\frac{\partial \delta \bm{x}_{\mathcal{E}'}}{\partial \bm{x}_{\mathcal{E}}^c} = -\mathbf{I}_{6\times 6}
\end{equation}
Substituting Eqs. (\ref{eq:abslute GVE}) and (\ref{eq:partials dxe dxed}) into Eq. (\ref{eq:roe chain rule}) and assuming zero interspacecraft separation yields
\begin{align}
\frac{d\delta \bm{x}_{\mathcal{E}'}}{dt} &= \mathbf{G}(\bm{x}_{\mathcal{E}}^d)  \underset{c \rightarrow d}{\mathbf R}\bm{d}_d  -  \mathbf{G}(\bm{x}_{\mathcal{E}}^c)\bm{d}_c\\
	&= \mathbf{G}(\bm{x}_{\mathcal{E}}^c)\delta\bm{d}
\end{align}
where $\delta\bm{d} = \bm{d}_d-\bm{d}_c$ is the difference between the deputy and chief perturbing accelerations expressed in the chief RTN frame. The matrix \mbox{\tiny$\raisebox{1.8pt}{$\underset{c\rightarrow d}{}$}$  \hspace{-.43cm}\small $\raisebox{2pt}{$\mathbf R$}$ \hspace{.01cm} \normalsize $\in \mathbb{R}^{3\times 3}$} rotates vectors from the chief RTN frame to the deputy RTN frame and is the identity matrix for zero interspacecraft separation. Thus the process noise mapping matrix is $\mathbf{\Gamma}_{\delta\mathcal{E}'} = \mathbf{G}(\bm{x}_{\mathcal{E}}^c)$, which matches that of the absolute equinoctial orbital elements shown in Eq. (\ref{eq:Gamma equinoctial}).

Since the state transition matrix and process noise mapping matrix of $\delta \bm{x}_{\mathcal{E}'}$ match that of the absolute equinoctial elements, the expressions for $\mathbf{Q}_{\delta\mathcal{E}',k}$ and $\mathbf{Q}_{\mathcal{E},k}$ are equivalent except that the expression for $\mathbf{Q}_{\delta\mathcal{E}',k}$ considers differential unmodeled accelerations. Thus, $\mathbf{Q}_{\delta\mathcal{E}',k}$ can be analytically approximated through Eq. (\ref{eq:Q split RTN}) and Eqs. (\ref{eq:QR equinoctial}-\ref{eq:Q equinoctial integration}) for a near-circular chief, Eqs. (\ref{eq:Q split RTN}) and (\ref{eq:Q equi small interval}) for a small propagation interval, or Eq. (\ref{eq:Q pert eccen}) by making the substitutions $\mathbf{Q}_{\mathcal{E},k} = \mathbf{Q}_{\delta\mathcal{E}',k}$ and $\widetilde{\mathbf{Q}}^\mathcal{R}= \widetilde{\mathbf{Q}}_{\delta}^\mathcal{R}$. Then following the approach in Eqs. (\ref{eq:P f(QR) abs 1}-\ref{eq:QI from Qpsi}), the process noise covariance of $\delta \bm{x}_{\mathcal{E}}$ is
\begin{equation}\label{eq:Qde small separation}
	\mathbf{Q}_{\delta\mathcal{E},k}  =  \mathbf{J}_{\delta\mathcal{E}}(t_k) \mathbf{Q}_{\delta\mathcal{E}',k} \mathbf{J}_{\delta\mathcal{E}}(t_k)^T
\end{equation}
where
\begin{equation}\label{eq:partials dxE wrt dxE'}
	\mathbf{J}_{\delta\mathcal{E}}(t_k) = \frac{\partial \delta\bm{x}_{\mathcal{E}}}{\partial \delta\bm{x}_{\mathcal{E}'}}\bigg\rvert_{\delta\bm{x}_{\mathcal{E}'}=\bm{0}}
		=
		\begin{bmatrix}
			\frac{1}{a_c} 					&\mathbf{0}_{1\times 4} 		&0\\
			0 										&\mathbf{0}_{1\times 4}		&1\\
			\mathbf{0}_{4\times 1}	&\mathbf{I}_{4\times 4} 		&\mathbf{0}_{4\times 1}
		\end{bmatrix}
\end{equation}

\section{Relative Process Noise Covariance for Large Separations}
If the interspacecraft separation is large, the linearization about zero interspacecraft separation utilized in the previous section may introduce significant errors. Instead of assuming small separations, this section constructs the process noise covariance of a relative spacecraft state by assuming that the interspacecraft separation is large enough that the unmodeled accelerations of the chief are uncorrelated with those of the deputy. This framework is flexible in that the relative state process noise covariance can be formed using any absolute state process noise covariance model, not just SNC models.

In general, 
\begin{align}\label{eq:rel state fun of d and c}
	\delta \bm{x}_{\alpha}(t_k) = \bm{f}(\bm{x}_{\beta}^d(t_k),\bm{x}_{\beta}^c(t_k))
\end{align}
where $\bm{f}(\cdot)$ is some function. The absolute states of the deputy and chief are $\bm{x}_{\beta}^d$ and $\bm{x}_{\beta}^c$ respectively where $\beta$ denotes the specific absolute state representation. Eq. (\ref{eq:rel state fun of d and c}) can be approximated through a first-order Taylor series expansion as
\begin{align}\label{eq:rel state taylor series}
	\delta \bm{x}_{\alpha}(t_k) = \delta \hat{\bm{x}}_{\alpha,k|k-1} + \mathbf{J}_d(t_k)(\bm{x}_{\beta}^d(t_k) - \hat{\bm{x}}_{\beta,k|k-1}^d ) + \mathbf{J}_c(t_k)(\bm{x}_{\beta}^c(t_k) - \hat{\bm{x}}_{\beta,k|k-1}^c ) 
\end{align}
which is valid for small deviations of the estimated chief and deputy absolute states from their corresponding true states. Here the partial derivatives
\begin{equation}
\begin{aligned}
\mathbf{J}_d(t_k) &= \frac{\partial\delta \bm{x}_{\alpha}(t_k)}{\partial \bm{x}_{\beta}^d(t_k)},
	\hspace{.2cm}&  
	 \mathbf{J}_c(t_k) &= \frac{\partial\delta \bm{x}_{\alpha}(t_k)}{\partial \bm{x}_{\beta}^c(t_k)}
\end{aligned}
\end{equation}
are evaluated at the mean state estimate. Using the linearization in Eq. (\ref{eq:rel state taylor series}), the propagated error covariance of $\delta \bm{x}_{\alpha}$ can be written as
\begin{align}
	\mathbf{P}_{\delta \alpha,k|k-1} &= \text{E}[(\delta\bm{x}_{\alpha}(t_{k}) - \delta\hat{\bm{x}}_{\alpha,k|k-1})(\delta\bm{x}_{\alpha}(t_{k}) - \delta\hat{\bm{x}}_{\alpha,k|k-1})^T]\\
	&= \mathbf{\Phi}_{\delta \alpha,k}\mathbf{P}_{\delta \alpha,k-1|k-1}\mathbf{\Phi}^T_{\delta \alpha,k} + \mathbf{Q}_{\delta \alpha,k}\label{eq:relative process noise covariance full}
\end{align}
where 
\begin{align}
	\mathbf{Q}_{\delta \alpha,k} &= \mathbf{J}_d (t_k) \mathbf{Q}_{\alpha,k}^d \mathbf{J}_d (t_k)^T + \mathbf{J}_c(t_k) \mathbf{Q}_{\alpha,k}^c \mathbf{J}_c (t_k)^T + \mathbf{J}_c(t_k)  \text{E}[\bm{w}_{\beta,k}^c\bm{w}_{\beta,k}^{d^T}] \mathbf{J}_d(t_k)^T\nonumber\\
	& \hspace{0.4cm} + (\mathbf{J}_c(t_k)  \text{E}[\bm{w}_{\beta,k}^c\bm{w}_{\beta,k}^{d^T}] \mathbf{J}_d(t_k)^T)^T \label{eq:relative Q full}
\end{align}
is process noise covariance of $\delta \bm{x}_{\alpha}$ at time step $k$. The discrete-time process noise of the absolute chief and deputy states are $\bm{w}_{\beta,k}^c$ and $\bm{w}_{\beta,k}^d$ respectively with associated covariances $\mathbf{Q}_{\beta,k}^c$ and $\mathbf{Q}_{\beta,k}^d$. In the particular case of SNC,
\begin{equation}\label{eq:relative Q SNC cross term}
	\text{E}[\bm{w}_{\beta,k}^c\bm{w}_{\beta,k}^{d^T}] = \int_{t_{k-1}}^{t_k} \mathbf{\Phi}_{\beta}^c(t_k,\tau)  \mathbf{\Gamma}_{\beta}^c(\tau)  \mathbf{\widetilde{Q}}_{cd} \mathbf{\Gamma}_{\beta}^d(\tau)^T  \mathbf{\Phi}_{\beta}^d(t_k,\tau)^T  d\tau
\end{equation}
where the superscripts $c$ and $d$ indicate matrices corresponding to the chief and deputy states respectively. Recall that $\mathbf{\widetilde{Q}}_{cd}$ is the cross power spectral density of the chief and deputy unmodeled accelerations where $\text{E}[\bm{\epsilon}_c(t) \bm{\epsilon}_d(\tau)^T] = \mathbf{\widetilde{Q}}_{cd}\delta(t-\tau)$. Note that $\bm{\epsilon}_c$ and $\bm{\epsilon}_d$ are expressed in the RTN frames of the chief and deputy respectively. If the interspacecraft separation is large, it can reasonably be assumed that the chief unmodeled accelerations are uncorrelated with those of the deputy such that $\text{E}[\bm{w}_{\beta,k}^c\bm{w}_{\beta,k}^{d^T}] = \bm{0}$. Using this assumption, Eq. (\ref{eq:relative Q full}) simplifies to
\begin{equation}\label{eq:relative Q large sep}
	\mathbf{Q}_{\delta \alpha,k} = \mathbf{J}_d (t_k) \mathbf{Q}_{\beta,k}^d \mathbf{J}_d (t_k)^T + \mathbf{J}_c (t_k) \mathbf{Q}_{\beta,k}^c \mathbf{J}_c (t_k)^T
\end{equation}
Eq. (\ref{eq:relative Q large sep}) can be applied regardless of how the process noise of the chief and deputy states are modeled provided that $\text{E}[\bm{w}_{\beta,k}^c\bm{w}_{\beta,k}^{d^T}]$ is small. Thus, $\mathbf{Q}_{\beta,k}^d$ and $\mathbf{Q}_{\beta,k}^c$ can be computed using any process noise covariance models for absolute spacecraft states, not just SNC models. Furthermore, Eq. (\ref{eq:relative Q large sep}) guarantees $\mathbf{Q}_{\delta \alpha,k}$ is positive semi-definite when $\mathbf{Q}_{\beta,k}^d$ and $\mathbf{Q}_{\beta,k}^c$ are positive semi-definite.


As an example, consider the relative Cartesian state $\delta \bm{x}_{\mathcal{R}}$ defined in Eq. (\ref{eq:relative cartesian state}), which is a function of the absolute, inertial Cartesian states of the chief and deputy denoted $\bm{x}_{\mathcal{I}}^c$ and $\bm{x}_{\mathcal{I}}^d$ respectively. The process noise covariance of $\delta \bm{x}_{\mathcal{R}}$ can be modeled through Eq. (\ref{eq:relative Q large sep}) where the process noise covariances of the chief and deputy inertial Cartesian states, $\mathbf{Q}_{\mathcal{I},k}^c$ and $\mathbf{Q}_{\mathcal{I},k}^d$, are separately computed using any absolute process noise covariance model such as Eq. (\ref{eq:Q kinematic inertial}) for a small propagation interval, Eqs. (\ref{eq:Q split RTN}) and (\ref{eq:Q HCW components rel rad}-\ref{eq:HCW integrals relative}) for a near-circular orbit, or Eq. (\ref{eq:Q pert eccen}). The required partial derivatives in Eq. (\ref{eq:relative Q large sep}), $\mathbf{J}_d(t_k) = \frac{\partial\delta \bm{x}_{\mathcal{R}}(t_k)}{\partial \bm{x}_{\mathcal{I}}^d(t_k)}$, are equivalent to Eq. (\ref{eq:partial dxR wrt dxI}). Analytical partial derivatives $\mathbf{J}_c(t_k) = \frac{\partial\delta \bm{x}_{\mathcal{R}}(t_k)}{\partial \bm{x}_{\mathcal{I}}^c(t_k)}$ are provided in Appendix B.



As another example, the process noise covariance of the relative orbital element state $\delta \bm{x}_{\mathcal{E}}$ defined in Eq. (\ref{eq:roe equinoctial}) can be modeled by Eq. (\ref{eq:relative Q large sep}) where $\mathbf{Q}_{\mathcal{E},k}^c$ and $\mathbf{Q}_{\mathcal{E},k}^d$ are each computed through any absolute process noise covariance model such as Eqs. (\ref{eq:Q split RTN}) and (\ref{eq:QR equinoctial}-\ref{eq:Q equinoctial integration}) for a near-circular orbit, Eqs. (\ref{eq:Q split RTN}) and (\ref{eq:Q equi small interval}) for a small propagation interval, or Eq. (\ref{eq:Q pert eccen}). In this case, the required partial derivatives in Eq. (\ref{eq:relative Q large sep}) are
\begin{equation}
\mathbf{J}_c(t_k) = \frac{\partial\delta \bm{x}_{\mathcal{E}}(t_k)}{\partial \bm{x}_{\mathcal{E}}^c(t_k)} = 
	\begin{bmatrix}
			\frac{-a_d}{a_c^2} 			&\mathbf{0}_{1\times 4} 		&0\\
			0 										&\mathbf{0}_{1\times 4}		&-1\\
			\mathbf{0}_{4\times 1}	&-\mathbf{I}_{4\times 4} 		&\mathbf{0}_{4\times 1}
		\end{bmatrix}
\end{equation}
and $\mathbf{J}_d(t_k) = \frac{\partial\delta \bm{x}_{\mathcal{E}}(t_k)}{\partial \bm{x}_{\mathcal{E}}^d(t_k)}$, which is equivalent to Eq. (\ref{eq:partials dxE wrt dxE'}).

\section{Numerical Validation}
This section validates the developed analytical process noise covariance models by comparing them to numerical solutions of the integrals they approximate for Earth-orbiting spacecraft trajectories. The truth orbit numerical integration scheme and perturbation models are summarized in Table \ref{tab:truth orbit info}. These perturbation models and integration scheme are also used to numerically evaluate the process noise covariance as well as the required state transition matrices using the relations
\begin{equation}
	\mathbf{\Phi}_{\alpha}(t_0,t_0) = \mathbf{I}, \hspace{.4cm} \dot{\mathbf{\Phi}}_{\alpha}(t,t_0) = \mathbf{A}_{\alpha}(t)\mathbf{\Phi}_{\alpha}(t,t_0)
\end{equation}
These numerical solutions of the process noise covariance are considered the reference truth. The simulated spacecraft mass is 100 kg, cross-section area is 1 m$^2$, drag coefficient is 1, and radiation pressure coefficient is 1.2. Results are also provided for a Keplerian reference truth where the perturbations in Table \ref{tab:truth orbit info} are omitted. The presented results for a Keplerian reference truth hold for any orbit semi-major axis because the time intervals associated with the results are normalized by the orbit period.

\vspace{0.8cm}
\begin{table}[!h]
\centering
\begin{threeparttable}
\caption{Truth orbit propagation parameters.}
\begin{tabular}{l l c c c c}
\toprule
\toprule
Parameter 						&Value\\
\toprule
Integration Scheme 			&Fourth-order Runge-Kutta\\
Integration Step Size 		&Fixed, 10 s\\
Initial Epoch						&J2000, 2000 January 1 12 hrs\\ 
Earth Gravity					&GGM05S\cite{ries_development_2016}, degree and order 30 \\
Atmospheric Density		&Harris-Priester\cite{montenbruck_satellite_force_2000}, considers atmosphere diurnal bulge \\
Third Body Gravity			&Sun and Moon point masses\\ 		
Solar Radiation Pressure	&Spacecraft cross-section normal to sun, no eclipses\\		
\bottomrule
\bottomrule
\end{tabular}
\label{tab:truth orbit info}
\end{threeparttable}
\end{table}
\vspace{0.4cm}

The accuracy of the analytical process noise covariance models is quantified using two error metrics. The first error metric is
\begin{equation}
	\Delta t_{min} = \text{min}\{\Delta t \ : \ |(Q_{ii}(\Delta t)^{\frac{1}{2}} - \hat{Q}_{ii}(\Delta t)^{\frac{1}{2}})/Q_{ii}(\Delta t)^{\frac{1}{2}}| \geq 0.1 \ \text{ for any } 1\leq i \leq 6\}
\end{equation}
which is the smallest propagation interval length for which the absolute fractional error in any modeled process noise standard deviation is greater than or equal to 0.1. The element in the $i$\textsuperscript{th} row and column of the process noise covariance as determined numerically and analytically are $Q_{ii}$ and $\hat{Q}_{ii}$ respectively. The second error metric is the maximum absolute fractional error in any modeled process noise standard deviation for any propagation interval length up to one orbit, which is defined as
\begin{equation}
	\delta_{max} = \text{max}\{|(Q_{ii}(\Delta t)^{\frac{1}{2}} - \hat{Q}_{ii}(\Delta t)^{\frac{1}{2}})/Q_{ii}(\Delta t)^{\frac{1}{2}}|  \ : 1\leq i \leq 6, \ 0\leq \Delta t \leq 1 \text{ orbit}\}
\end{equation}





\subsection{Absolute Process Noise Covariance}
The analytical process noise covariance models of the absolute Cartesian state $\bm{x}_{\mathcal{I}}$ and equinoctial element state $\bm{x}_{\mathcal{E}}$ are each compared against the corresponding numerical solution of Eq. (\ref{eq:absolute process noise covariance}). While the initial  osculating eccentricity is varied, the orbit inclination and right ascension of the ascending node are each $45\degree$. The initial mean longitude is $90\degree$. The initial periapsis radius is fixed at 6900 km in order to ensure the trajectory is feasible and does not pass through the Earth for any of the considered eccentricities. Thus, the semi-major axis and orbit period vary depending on the eccentricity. The small periapsis radius leads to large perturbing accelerations due to nonspherical gravity and atmospheric drag. The spacecraft trajectory nominally starts at periapsis, although periapsis is not well defined for small eccentricities, in order to maximize the effects of eccentricity. 
 
The absolute power spectral density is as shown in Eq. (\ref{eq:psd RTN}). In one scenario, $\widetilde{Q}_r$, $\widetilde{Q}_t$, and $\widetilde{Q}_n$ are all equal to a single value denoted $\widetilde{Q}^*$. In a second scenario, $\widetilde{Q}_r = \widetilde{Q}_n = \widetilde{Q}^*$ and $\widetilde{Q}_t = 10\widetilde{Q}^*$, which is representative of a case where atmospheric drag is the largest source of unmodeled accelerations. Notice in Eq. (\ref{eq:absolute process noise covariance}) that each element of the process noise covariance is a linear function of the elements of the power spectral density. Since the nonzero elements of the power spectral density are each a linear function of $\widetilde{Q}^*$, each process noise standard deviation is some scalar times $(\widetilde{Q}^{*})^{\frac{1}{2}}$ for both the numerical and considered analytical solutions. Thus the fractional error of each modeled process noise standard deviation is the same regardless of the simulated value of $\widetilde{Q}^*$. 

\begin{figure}[H]
\centering
\subfigure[Equal Noise in Each RTN Axis]{\label{fig:Delta t}
\includegraphics[width=.485\linewidth,trim=0 -5 -5 -10,clip]{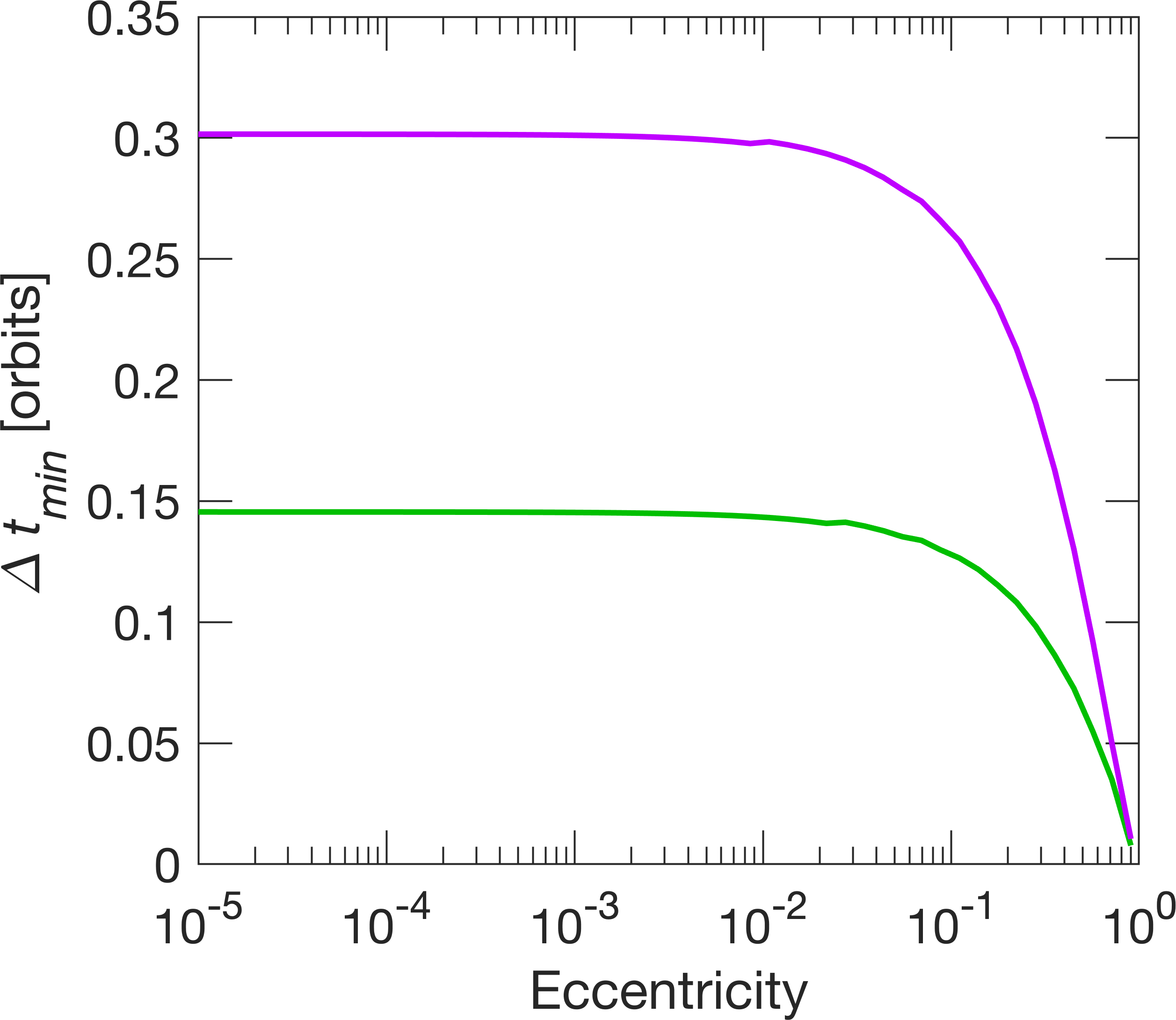}}
\subfigure[Transverse Dominant Noise]{\label{fig:Delta t Tdom}
\includegraphics[width=.485\linewidth,trim=-5 -5 0 -10,clip]{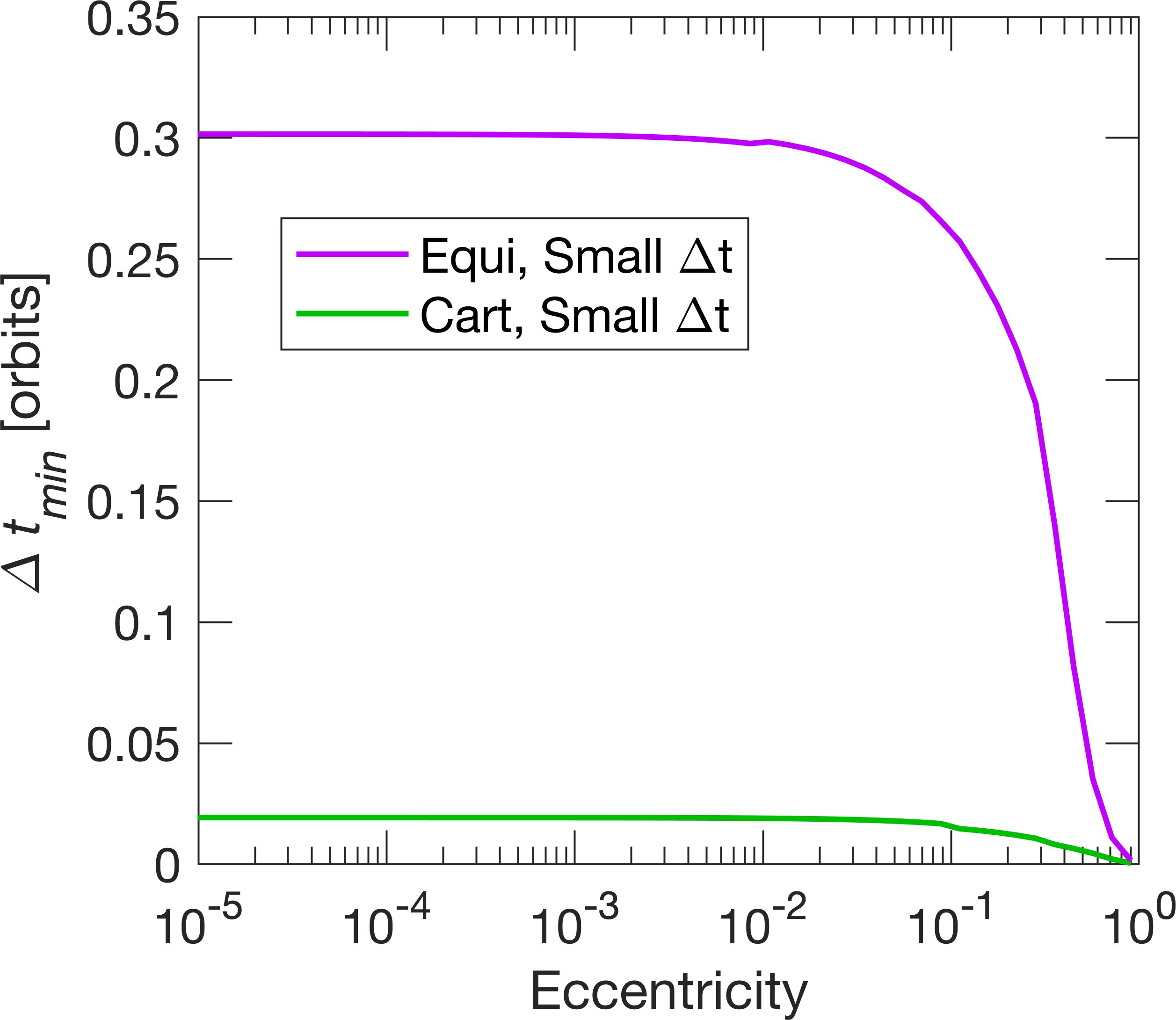}}
\caption{Comparison of $\Delta t_{min}$ for the new process noise covariance model in Eq. (\ref{eq:Q equi small interval}) for absolute equinoctial elements (purple) and the widely used kinematic model in Eq. (\ref{eq:Q kinematic inertial}) for an absolute, inertial Cartesian state (green).}
\label{fig:absolute tMin}
\end{figure}

The new model for equinoctial elements in Eq. (\ref{eq:Q equi small interval}) and the widely used model for an inertial Cartesian state in Eq. (\ref{eq:Q kinematic inertial}) both assume a small propagation interval $\Delta t_k$. To determine the range of propagation interval lengths for which these models can be accurately applied, $\Delta t_{min}$ is plotted in Figure \ref{fig:absolute tMin} as a function of eccentricity. The fractional error in the process noise standard deviation of the equinoctial element $h$ is as large as $0.14$ for $\Delta t < 0.13$ orbits. However, since the initial true longitude is $90\degree$, the process noise standard deviation of $h$ is very small for $\Delta t < 0.13$ orbits. Thus, the absolute error is very small, and the process noise standard deviation of $h$ for $\Delta t < 0.13$ orbits is neglected in Figure \ref{fig:absolute tMin}. The spacecraft angular velocity is greatest at periapsis and grows with increasing eccentricity, more quickly breaking the assumptions of the models in Eqs. (\ref{eq:Q equi small interval}) and (\ref{eq:Q kinematic inertial}). Thus, $\Delta t_{min}$ tends to decrease for larger eccentricities in Figure \ref{fig:absolute tMin} for both models, and both models generally perform better when the propagation interval is not near periapsis. The equinoctial element model outperforms the Cartesian model. Interestingly, $\Delta t_{min}$ for the Cartesian model decreases significantly when the noise is predominantly in the transverse axis. The results in Figure \ref{fig:absolute tMin} are nearly identical to when the reference truth is computed assuming Keplerian motion.

\vspace{.5cm}

\begin{figure}[H]
\centering
\subfigure[Keplerian Reference Truth]{\label{fig:absolute eMax kep}
\includegraphics[width=.485\linewidth,trim=0 -5 0 0,clip]{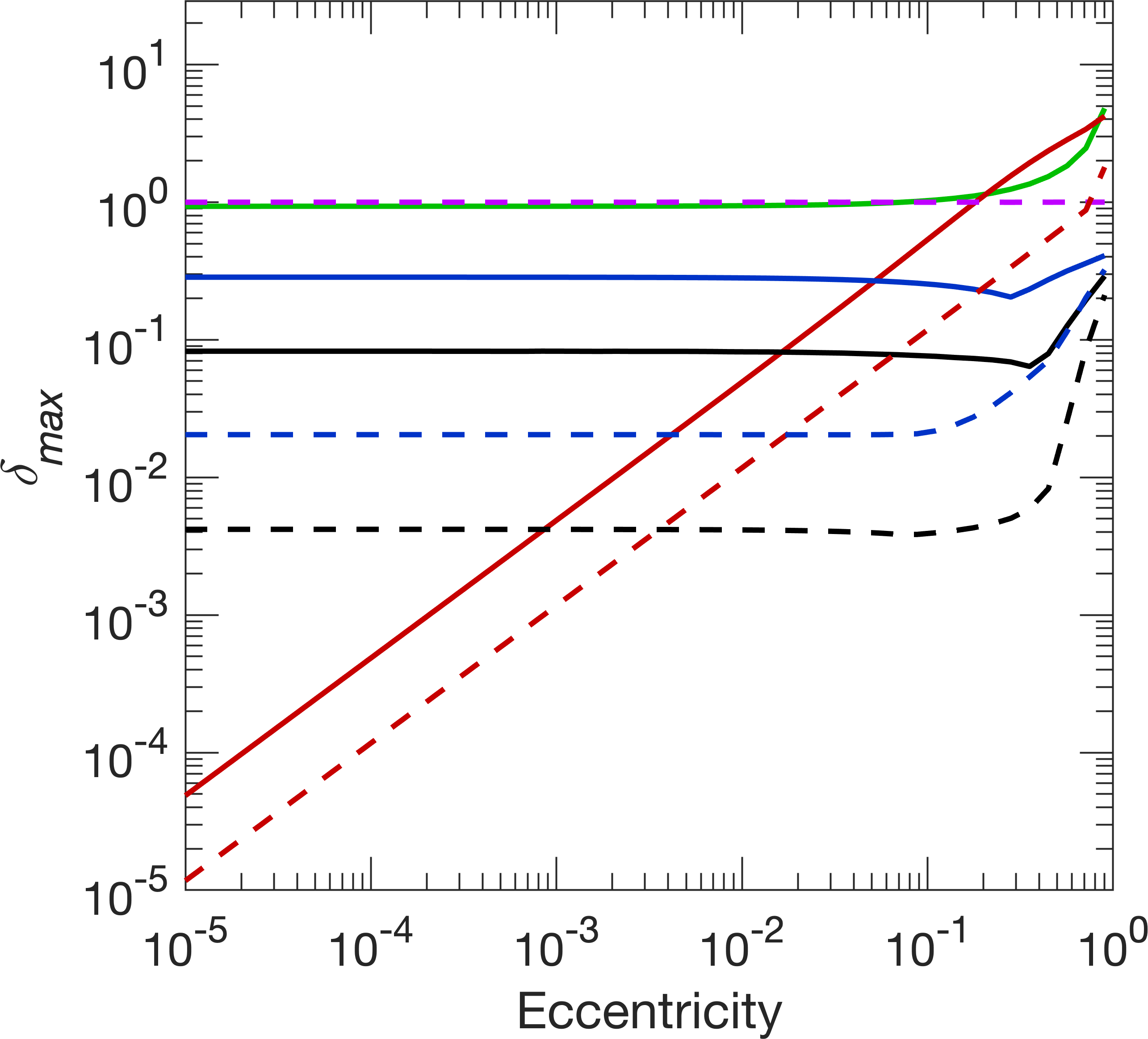}}
\subfigure[Perturbed Reference Truth (See Table \ref{tab:truth orbit info})]{\label{fig:absolute eMax pert}
\includegraphics[width=.485\linewidth,trim=0 -5 0 0,clip]{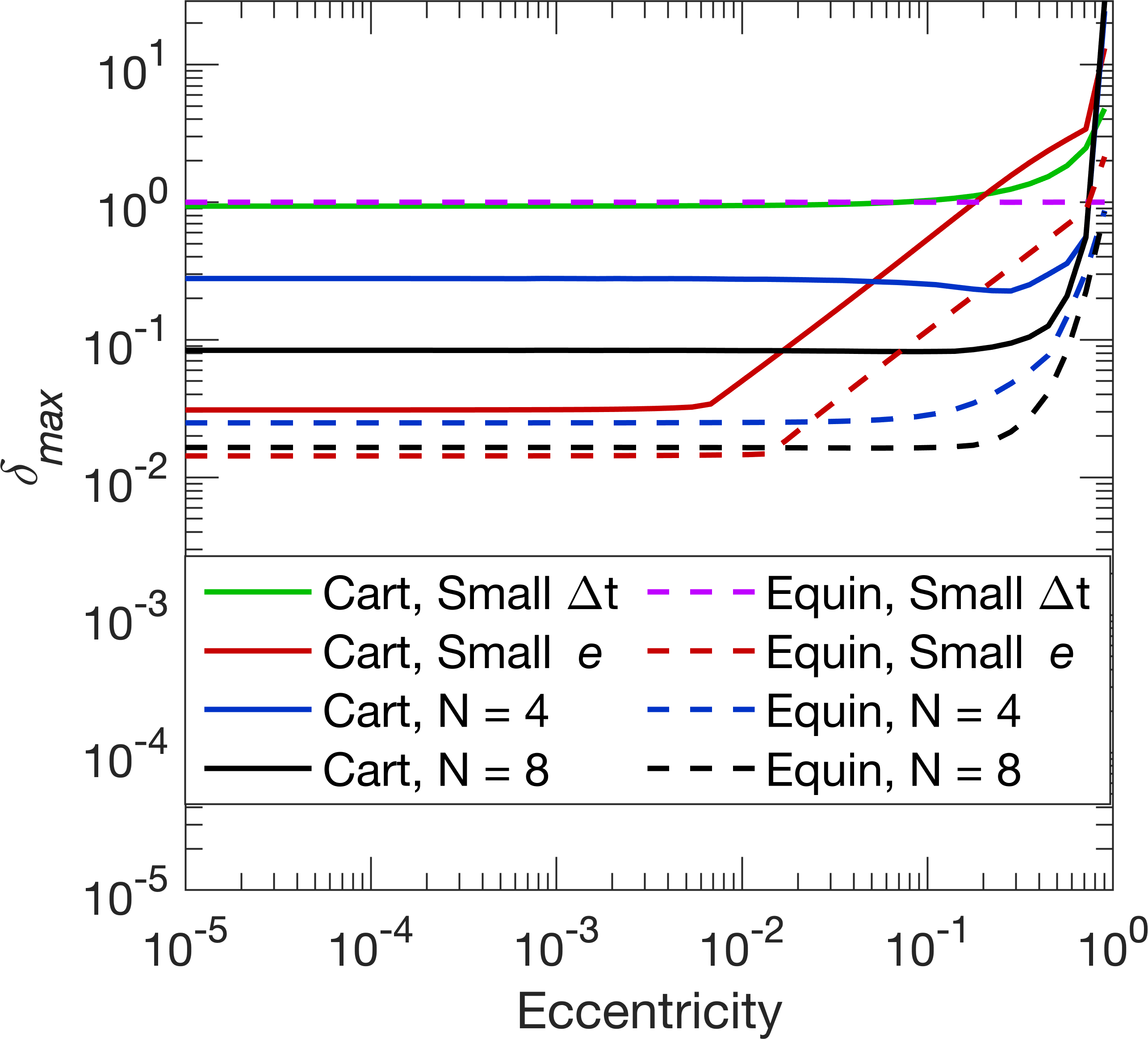}}
\caption{Comparison of $\delta_{max}$ for Eq. (\ref{eq:Q kinematic inertial}) (green) and Eq. (\ref{eq:Q equi small interval}) (purple) that assume a small propagation interval, Eq. (\ref{eq:QI from Qpsi}) (red solid) and Eqs. (\ref{eq:QR equinoctial}-\ref{eq:Q equinoctial integration}) (red dashed) that assume a circular orbit, and Eq. (\ref{eq:Q pert eccen}) for either four (blue) or eight (black) subintervals.}
\label{fig:absolute emax}
\end{figure}

The new analytical models for absolute equinoctial elements and Cartesian coordinates shown in Eqs. (\ref{eq:QR equinoctial}-\ref{eq:Q equinoctial integration}) and Eq. (\ref{eq:QI from Qpsi}) respectively assume a circular orbit. In order to quantify the maximum orbit eccentricity for which these models can be applied, $\delta_{max}$ is plotted in Figure \ref{fig:absolute emax} as a function of eccentricity for the case of equal noise in each axis. The results of the transverse dominant noise scenario were similar. The equinoctial model is valid for larger eccentricities than the Cartesian model. For example, $\delta_{max}$ is less than 0.1 for eccentricities less than about 0.09 and 0.02 for the equinoctial and Cartesian models respectively for both the Keplerian and perturbed reference truths. Since $\delta_{max}$ is the largest fractional error over a full orbit, the fractional errors are smaller on average. Depending on the particular application, $\delta_{max} \leq 0.1$ may provide sufficient accuracy. The presented model for equinoctial elements in Eq. (\ref{eq:Q equi small interval}) and the widely used kinematic model for an inertial Cartesian state in Eq. (\ref{eq:Q kinematic inertial}) are also included in Figure \ref{fig:absolute emax} for reference. Since both of these models assume a small propagation interval and $\delta_{max}$ considers intervals up to an entire orbit, it is not surprising that $\delta_{max}$ is large for these models.

Figure \ref{fig:absolute emax} also shows the performance of Eq. (\ref{eq:Q pert eccen}), which breaks the propagation interval into smaller subintervals, applied to absolute inertial Cartesian coordinates and equinoctial elements using either four or eight subintervals. For Cartesian coordinates, $\mathbf{Q}_{\alpha}(\bar{t}_i,\bar{t}_{i-1})$ is computed using Eq. (\ref{eq:Q kinematic inertial}), and $\mathbf{\Phi}_{\alpha}(t_k,\bar{t}_i)$ is modeled using the partial derivatives in Eqs. (\ref{eq:YA partials 1}-\ref{eq:YA partials 2}) and the Yamanaka-Ankerson state transition matrix formulated using argument of latitude\cite{willis_analytical_2021}. For equinoctial elements,  Eq. (\ref{eq:Q equi small interval}) is used to compute $\mathbf{Q}_{\alpha}(\bar{t}_i,\bar{t}_{i-1})$, and Eq. (\ref{eq:stm oe}) is used for $\mathbf{\Phi}_{\alpha}(t_k,\bar{t}_i)$. For simplicity, the subintervals are evenly spaced in time. However, further accuracy could likely be achieved with the same number of subintervals by concentrating more of the subintervals toward the beginning of the propagation interval since the earlier terms of the sum in Eq. (\ref{eq:Q pert eccen}) tend to be larger due to the propagation of $\mathbf{Q}_{\alpha}(\bar{t}_i,\bar{t}_{i-1})$ over a larger interval $[\bar{t}_i,t_k]$. For eccentric orbits, it may also be beneficial to concentrate more of the subintervals near periapsis since that is where the spacecraft angular velocity is greatest, more quickly violating the assumptions of Eqs. (\ref{eq:Q kinematic inertial}) and (\ref{eq:Q equi small interval}). As a worst case scenario, the states at each time $\bar{t}_i$, which are required to compute each $\mathbf{\Phi}_{\alpha}(t_k,\bar{t}_i)$ and $\mathbf{Q}_{\alpha}(\bar{t}_i,\bar{t}_{i-1})$ in Eq. (\ref{eq:Q pert eccen}), are computed through a Keplerian propagation of the state from the beginning of the propagation interval. Figure \ref{fig:absolute emax} shows that increasing the number of subintervals $N$ increases accuracy. The equinoctial elements model is more accurate than the Cartesian model given the same number of subintervals because Eq. (\ref{eq:Q equi small interval}) is accurate over longer intervals than Eq. (\ref{eq:Q kinematic inertial}) as shown in Figure \ref{fig:absolute tMin}. Even for the considered highly perturbed orbits, the models that use Eq. (\ref{eq:Q pert eccen}) and neglect perturbations provide sufficient accuracy for many applications for propagation intervals of an orbit or less. The accuracy of these models can be further improved by using $\mathbf{\Phi}_{\alpha}(t_k,\bar{t}_i)$ that incorporate perturbations in Eq. (\ref{eq:Q pert eccen}), which is important for propagation intervals of multiple orbits.

\vspace{.5cm}

\subsection{Relative Process Noise Covariance}
Consider two spacecraft where the chief has an initial osculating eccentricity of $10^{-3}$ and an inclination and right ascension of the ascending node of $45\degree$. The initial chief osculating equinoctial elements are $[a,f,g,h,k,\lambda] = [6900 \ \text{km},0,10^{-3},0.2929,0.2929,135\degree]$. The initial deputy mean equinoctial elements all match those of the chief with the exception of the mean longitude, which is varied. As a result, the two spacecraft are primarily separated in the along-track direction as described by the initial mean relative mean longitude, $\overline{\delta \lambda}_0$, and the secular drift in $\delta \lambda$ is minimal. The transformation between mean and osculating orbital elements is done using the first-order $J_2$ mapping presented by Schaub and Junkins\cite{schaub_analytical_2014}. The unmodeled acceleration power spectral density for each spacecraft is as shown in Eq. (\ref{eq:psd RTN}) where each diagonal element is equal to $\widetilde{Q}^*$ for both spacecraft. The presented results hold for any $\widetilde{Q}^*$. In reality, the degree of correlation between the unmodeled accelerations of two spacecraft depends on the fidelity of the modeled dynamics as well as the orbit geometry. Here, the cross power spectral density of the chief and deputy unmodeled accelerations is modeled as $\widetilde{\mathbf{Q}}_{cd} = (0.9e^{-\overline{\delta \lambda}_0/10^{-3}}) \widetilde{\mathbf{Q}}_{c}$. As a result, $\widetilde{\mathbf{Q}}_{cd}$ decays exponentially as $\overline{\delta \lambda}_0$ increases. The rate of decay is dictated by the constant $10^{-3}$, which was chosen for illustrative purposes. The constant 0.9 models the fact that the unmodeled accelerations of the two spacecraft are not perfectly correlated near zero separation due to small differences in attitude and spacecraft physical properties. The reference truth process noise covariance of the relative states $\delta\bm{x}_\mathcal{R}$ and $\delta\bm{x}_\mathcal{E}$ are computed through Eq. (\ref{eq:relative Q full}). For $\delta\bm{x}_\mathcal{R}$, the partial derivatives $\frac{\partial\delta \bm{x}_{\mathcal{R}}(t_k)}{\partial \bm{x}_{\mathcal{I}}^c(t_k)}$ in Eq. (\ref{eq:relative Q full}) are approximated through central finite differences using the exact angular velocity vector of the RTN frame with respect to the inertial frame expressed in the RTN frame $\bm{w} = [r_c(\ddot{\bm{r}}_c\cdot\hat{\bm{n}})/L \ \ 0 \ \ L/r_c^2]^T$, which holds for arbitrary spacecraft motion\cite{casotto_equations_2016}. The reference truth process noise covariances of the absolute chief and deputy states as well as their cross covariance, which are employed in Eq. (\ref{eq:relative Q full}), are obtained through numerical integration of Eqs. (\ref{eq:absolute process noise covariance}) and (\ref{eq:relative Q SNC cross term}) according to Table \ref{tab:truth orbit info} both with and without perturbations. 

Two approaches were taken to develop analytical process noise covariance models for relative spacecraft states by assuming either small or large interspacecraft separation. For the process noise covariance of $\delta\bm{x}_\mathcal{R}$, the absolute model in Eq. ({\ref{eq:QI from Qpsi}) is applied in Eq. (\ref{eq:Qdr small separation}) for the small separation framework and Eq. (\ref{eq:relative Q large sep}) for the large separation framework. The analytical partial derivatives $\frac{\partial\delta \bm{x}_{\mathcal{R}}(t_k)}{\partial \bm{x}_{\mathcal{I}}^c(t_k)}$ in Appendix B are used in Eq. (\ref{eq:relative Q large sep}). For the process noise covariance of $\delta\bm{x}_\mathcal{E}$, the absolute model in Eqs. (\ref{eq:QR equinoctial}-\ref{eq:Q equinoctial integration}) is applied in Eq. (\ref{eq:Qde small separation}) for the small separation framework and Eq. (\ref{eq:relative Q large sep}) for the large separation framework. Since the employed absolute models in Eq. ({\ref{eq:QI from Qpsi}) and Eqs. (\ref{eq:QR equinoctial}-\ref{eq:Q equinoctial integration}) assume circular orbits, the modeled relative process noise covariance incurs some error due to eccentricity whether the small or large separation approach is applied. The large separation framework errors increase for small separations due to the assumption that the chief unmodeled accelerations are uncorrelated with those of the deputy. On the other hand, the small separation framework fully accounts for the correlation between the chief and deputy unmodeled accelerations as shown in Eq. (\ref{eq:cross power spectral density}). However, the small separation framework errors grow as interspacecraft separation increases due to the linearization about zero interspacecraft separation.

In order to quantify these errors, $\delta_{max}$ is plotted in Figure \ref{fig:relative emax} for the considered analytical models. These results demonstrate the inherent trade-off in the large and small separation frameworks, which are each valid for a different range of interspacecraft separations. Most notably, the small separation framework can be used for much larger interspacecraft separations for equinoctial elements than for Cartesian coordinates. Specifically, $\delta_{max} \leq 0.1$ for the Cartesian and equinoctial small separation models for $\delta\lambda \leq 2\times 10^{-3}$ and $\delta\lambda \leq 1\times 10^{-1}$ respectively for both the Keplerian and perturbed reference truths. These angles can be transformed to approximate interspacecraft distances of 13.8 km and 690 km respectively through multiplication with the orbit semi-major axis. Even for the perturbed reference truth, $\delta_{max} < 0.013$ can be achieved for any $\overline{\delta \lambda}_0$ for an equinoctial element state by using either the small or large separation framework. On the other hand, the lowest $\delta_{max}$ that can be achieved for a Cartesian state using the small and large separation frameworks is as high as 0.1 for the considered $\overline{\delta \lambda}_0$.

\begin{figure}[H]
\centering
\subfigure[Keplerian Reference Truth]{\label{fig:relative eMax kep}
\includegraphics[width=.485\linewidth,trim=0 -5 0 0,clip]{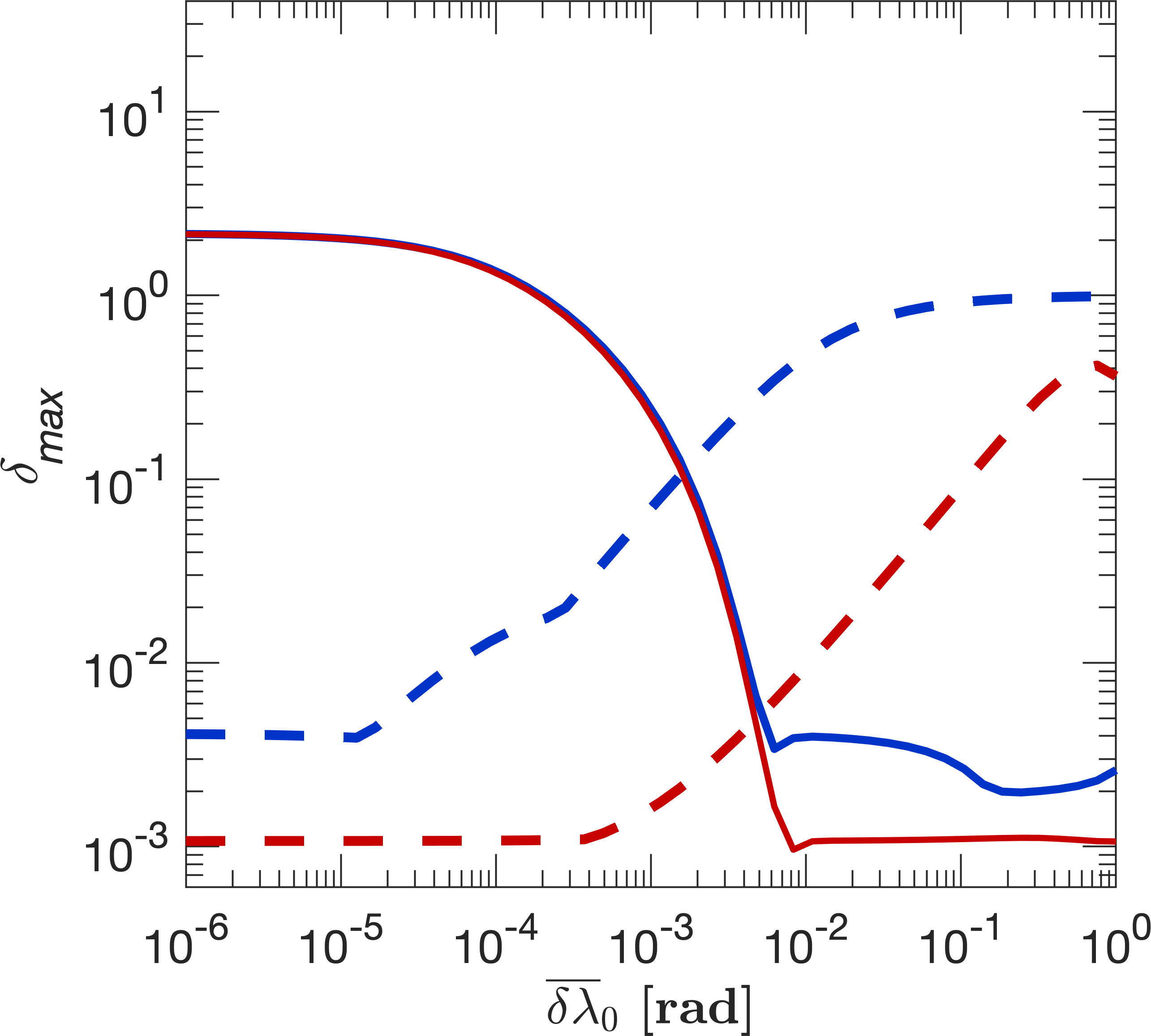}}
\subfigure[Perturbed Reference Truth (See Table \ref{tab:truth orbit info})]{\label{fig:relative eMax pert}
\includegraphics[width=.485\linewidth,trim=0 -5 0 0,clip]{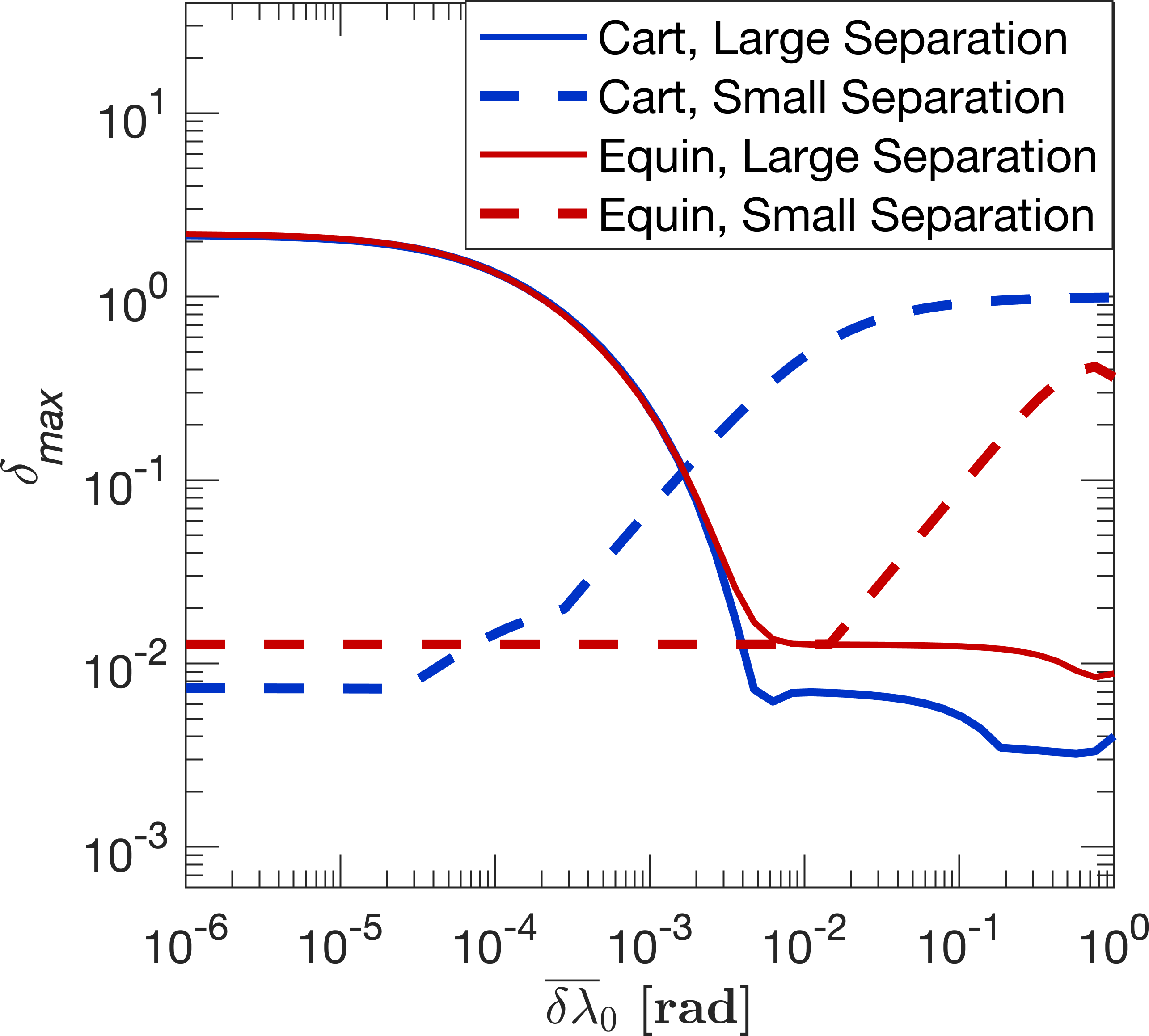}}
\caption{Comparison of $\delta_{max}$ when assuming small separation (dashed) and large separation (solid) for relative Cartesian (blue) and equinoctial element (red) states.}
\label{fig:relative emax}
\end{figure}
\vspace{.5cm}

\section{Conclusions}
Accurate process noise modeling is essential for optimal orbit determination in a discrete-time Kalman filtering framework and can improve satellite conjunction analysis. A common approach to process noise modeling called state noise compensation (SNC) treats the process noise as zero-mean Gaussian white noise unmodeled accelerations. The resulting process noise covariance can be evaluated numerically. However, analytical solutions are desirable because they are computationally efficient and provide insight into system behavior. One analytical SNC model is widely used that assumes kinematic motion, but it is restricted to an absolute Cartesian state and small propagation intervals. Furthermore, analytical SNC models for absolute orbital element states are not currently available and there has been little work on process noise covariance modeling for relative spacecraft states. This paper addresses these gaps in the state of the art.

First, a new analytical SNC process noise covariance model is developed for an absolute, inertial Cartesian spacecraft state by leveraging the well-known solution of the Hill-Clohessy-Wiltshire equations. For small orbit eccentricities, this new model can be used for significantly longer propagation intervals than the widely used kinematic model. An approach is presented to obtain analytical SNC process noise covariance models for absolute orbital element states for small eccentricities, and then a second approach is presented that is valid for eccentric orbits over small propagation intervals. The developed absolute Cartesian and orbital element models are then extended to long propagation intervals for perturbed, eccentric orbits by leveraging spacecraft relative dynamics models. Subsequently, two frameworks are established for modeling the process noise covariance of relative spacecraft states by assuming either small or large interspacecraft separations. The large interspacecraft separation framework is particularly flexible in that it can leverage any absolute state process noise covariance model, not just SNC models such as those developed in this paper. All the presented process noise covariance models are analytical, making them much more computationally efficient than numerical solutions, and are guaranteed positive semi-definite. Overall, the orbital element models tend to outperform the Cartesian models. Nevertheless, the newly developed analytical process noise covariance models for absolute and relative spaceraft states provide sufficient accuracy for many applications for both Cartesian and orbital element state representations.

\section*{Appendix A: Gauss Variational Equations}
\setcounter{equation}{0}
\renewcommand\theequation{A.\arabic{equation}}
The Gauss Variational Equations formulated in equinoctial elements are provided here for reference\cite{battin_introduction_1987,roth_gaussian_1985}. The parameters in Eq. (\ref{eq:Gamma equinoctial}) are
\begin{alignat}{3}
\bar{A} &= \frac{2a^2}{L}(f\sin l - g\cos l)			 							&&\bar{B} = \frac{2a^2W}{L}					&&\bar{C} = \sqrt{\frac{p}{\mu}}\sin l \nonumber\\
\bar{D} &= \sqrt{\frac{p}{\mu}}\frac{f+(1+W)\cos l}{W}	&&\bar{E} = \sqrt{\frac{p}{\mu}}\frac{g(k\cos l -h\sin l)}{W}							&&\bar{F} = -\sqrt{\frac{p}{\mu}}\cos l\nonumber\\
\bar{G} &= \sqrt{\frac{p}{\mu}}\frac{g+(1+W)\sin l}{W}\quad  &&\bar{H} = \sqrt{\frac{p}{\mu}}\frac{f(h\sin l-k\cos l)}{W}  \quad	&&\bar{I} = \sqrt{\frac{p}{\mu}}\frac{(1+h^2+k^2)\cos l}{2W}\nonumber\\
\bar{J} &= \sqrt{\frac{p}{\mu}}\frac{(1+h^2+k^2)\sin l}{2W}	&&\bar{K} = -\sqrt{\frac{p}{\mu}}\left(\frac{W-1}{1+\eta} + \frac{2\eta}{W} \right)						&&\bar{L} = -\frac{L(1+W)(g\cos l-f\sin l)}{\mu W(1+\eta)} \nonumber\\
\bar{M} &=-\frac{L(k\cos l-h\sin l)}{\mu W} 			&&W = 1+f\cos l+g\sin l				&& \eta = \sqrt{1-e^2} \label{eq:GVE equinoctial}
\end{alignat}
Note that Roth\cite{roth_gaussian_1985} has a sign error on $\bar{H}$.

\section*{Appendix B: Additional Partial Derivatives}
\setcounter{equation}{0}
\renewcommand\theequation{B.\arabic{equation}}
This appendix provides analytical partial derivatives 
\begin{equation}\label{eq:appendix partials}
\frac{\partial\delta \bm{x}_{\mathcal{R}}(t_k)}{\partial \bm{x}_{\mathcal{I}}^c(t_k)} = 
\begin{bmatrix}
	\frac{\partial \delta \bm{r}}{\partial \bm{r}_c} 	&\frac{\partial \delta \bm{r}}{\partial \bm{v}_c}\\
	\frac{\partial \delta \bm{v}}{\partial \bm{r}_c} 	&\frac{\partial \delta \bm{v}}{\partial \bm{v}_c}\\
\end{bmatrix}
\end{equation}
which were validated against a central finite difference approximation. First, note that the rotation matrix that maps vectors from the inertial frame to the RTN frame is
\begin{equation}
	\underset{\mathcal{I} \rightarrow \mathcal{R}}{\mathbf R} = 
	\begin{bmatrix}
		\hat{\bm{r}}^T\\
		\hat{\bm{t}}^T\\
		\hat{\bm{n}}^T\\
	\end{bmatrix}
	=
	\begin{bmatrix}
		\frac{\bm{r}_c^T}{||\bm{r}_c||}\\
		(\hat{\bm{n}} \times\hat{\bm{r}})^T\\
		\frac{(\bm{r}_c \times \bm{v}_c)^T}{||\bm{r}_c \times \bm{v}_c||}
	\end{bmatrix}
\end{equation}
where $\hat{\bm{r}}$, $\hat{\bm{t}}$, and $\hat{\bm{n}}$ are the unit vectors of the RTN frame expressed in the inertial frame. The angular velocity vector of the RTN frame with respect to the inertial frame expressed in the RTN frame is $\bm{w} = [r_c(\ddot{\bm{r}}_c\cdot\hat{\bm{n}})/L \ \ 0 \ \ L/r_c^2]^T$, which holds for arbitrary spacecraft motion\cite{casotto_equations_2016}. Neglecting out-of-plane accelerations, which are typically very small compared to the in-plane accelerations, the angular velocity vector is $\bm{w} = [0 \ \ 0 \ \ L/r_c^2]^T$. This approximation avoids any partial derivatives of the chief inertial acceleration $\ddot{\bm{r}}_c$. The chief specific angular momentum vector is denoted $\mathbf{L} = \bm{r}_c\times\bm{v}_c$, and its magnitude is $L = ||\bm{L}||$. Now the partial derivatives in Eq. (\ref{eq:appendix partials}) can be deduced from Eq. (\ref{eq:relative cartesian state}), which yields
\begin{align}
\frac{\partial \delta \bm{r}}{\partial \bm{r}_c}  &= 
	\begin{bmatrix}
		\frac{-\bm{r}_d^T\bm{r}_c\bm{r}_c^T}{r_c^3} + \frac{\delta\bm{r}_{\mathcal{I}}^T}{r_c}\\
		\bm{r}_d^T\mathbf{K}_1\\
		\frac{\bm{L}^T\bm{r}_d(\bm{L}\times \bm{v}_c)^T}{L^3} + \frac{(\bm{v}_c\times \bm{r}_d)^T}{L}
	\end{bmatrix}\label{eq:appendix partials 1}\\		
\frac{\partial \delta \bm{r}}{\partial \bm{v}_c}  &= 
	\begin{bmatrix}
		\mathbf{0}_{1\times 3}\\
		\delta\bm{r}_{\mathcal{I}}^T\mathbf{K}_2\\
		\delta\bm{r}_{\mathcal{I}}^T(\frac{\bm{r}_c^{\times}}{L} - \frac{\bm{L}(\bm{L}\times\bm{r}_c)^T}{L^3})
	\end{bmatrix}\label{eq:appendix partials 2}\\
\frac{\partial \delta \bm{v}}{\partial \bm{r}_c}  &=
	\begin{bmatrix}
		\frac{\delta\bm{v}_{\mathcal{I}}^T(r_c^2\mathbf{I}_{3\times3}-\bm{r}_c\bm{r}_c^T)}{r_c^3} - \delta\bm{r}_{\mathcal{I}}^T(\hat{\bm{t}}(\frac{(\bm{L}\times\bm{v}_c)^T}{Lr_c^2} + \frac{2L\bm{r}_c^T}{r_c^4}) - \frac{L\mathbf{K}_1}{r_c^2}) - \frac{L\hat{\bm{t}}^T}{r_c^2}\\
		\delta\bm{v}_{\mathcal{I}}^T\mathbf{K}_1 - \frac{L(\bm{r}_d-2\bm{r}_c)^T}{r_c^3} + \bm{r}_c^T\delta\bm{r}_{\mathcal{I}}(\frac{(\bm{L}\times \bm{v}_c)^T}{Lr_c^3} + \frac{3L\bm{r}_c^T}{r_c^5})\\
		\delta\bm{v}_{\mathcal{I}}^T(\frac{\bm{L}(\bm{L}\times \bm{v}_c)^T}{L^3} - \frac{\bm{v}_c^{\times}}{L})
	\end{bmatrix}\label{eq:appendix partials 3}\\
\frac{\partial \delta \bm{v}}{\partial \bm{v}_c}  &= 
	\begin{bmatrix}
		-\hat{\bm{r}}^T + \frac{\delta\bm{r}_{\mathcal{I}}^T}{r_c^2}(\frac{\hat{\bm{t}}(\bm{L}\times \bm{r}_c)^T}{L} + L\mathbf{K}_2)\\
		-\hat{\bm{t}}^T + \delta\bm{v}_{\mathcal{I}}^T\mathbf{K}_2 - \frac{\delta\bm{r}_{\mathcal{I}}^T\bm{r}_c(\bm{L}\times \bm{r}_c)^T}{Lr_c^3}\\
		\frac{-(\bm{L}^T+(\bm{r_c}\times\delta\bm{v}_{\mathcal{I}})^T)}{L} - \frac{\bm{L}^T\delta\bm{v}_{\mathcal{I}}(\bm{L}\times \bm{r}_c)^T}{L^3}\label{eq:appendix partials 4}
	\end{bmatrix}
\end{align}
Recall that the superscript $\times$ indicates a cross-product matrix as illustrated in Eq. (\ref{eq:w cross product mat}). The difference between the deputy and chief inertial position and velocity vectors are $\delta \bm{r}_{\mathcal{I}} = \bm{r}_d-\bm{r}_c$ and $\delta \bm{v}_{\mathcal{I}} = \bm{v}_d-\bm{v}_c$. The auxilliary variables  $\mathbf{K}_1, \mathbf{K}_2 \in \mathbb{R}^{3\times 3}$ are
\begin{align}
	\mathbf{K}_1 &= \frac{\partial \hat{\bm{t}}}{\partial \bm{r}_c} = \hat{\bm{n}}^\times \left(\frac{\mathbf{I}_{3\times 3}}{r_c} - \frac{\bm{r}_c\bm{r}_c^T}{r_c^3}\right) - \hat{\bm{r}}^\times \left(\frac{\bm{L}(\bm{L}\times\bm{v}_c)^T}{L^3} - \frac{\bm{v}_c^{\times}}{L}\right)\label{eq:gamma1}\\
	\mathbf{K}_2 &= \frac{\partial \hat{\bm{t}}}{\partial \bm{v}_c} = - \hat{\bm{r}}^\times \left(\frac{\bm{r}_c^{\times}}{L} - \frac{\bm{L}(\bm{L}\times\bm{r}_c)^T}{L^3}\right)\label{eq:gamma2}
\end{align}
In order to account for out-of-plane accelerations, the partial derivatives of $([\frac{r_c}{L}(\ddot{\bm{r}}_c\cdot\hat{\bm{n}}) \ \ 0 \ \ 0]^T)^{\times}\underset{\mathcal{I} \rightarrow \mathcal{R}}{\mathbf R}\delta\bm{r}_{\mathcal{I}}$ with respect to $\bm{r}_c$ and $\bm{v}_c$ would be subtracted from the right hand sides of Eqs. (\ref{eq:appendix partials 3}) and (\ref{eq:appendix partials 4}) respectively.

\section*{Acknowledgments}
This material is based upon work supported by the National Science Foundation (NSF) Graduate Research Fellowship Program under Grant No. DGE-1656518. Any opinions, findings, and conclusions or recommendations expressed in this material are those of the authors and do not necessarily reflect the views of the National Science Foundation. This work is also supported by the NASA Small Spacecraft Technology Program cooperative agreement number 80NSSC18M0058 for contributions to the Autonomous Nanosatellite Swarming (ANS) Using Radio-Frequency and Optical Navigation project. Additionally, the authors wish to thank the Achievement Rewards for College Scientists (ARCS) Foundation for their support.


\bibliographystyle{elsarticle-num}   
\bibliography{references}   

\end{document}